\documentclass{artjlt}

\usepackage{amsmath,amsfonts,amssymb}

\title{Structures of  Nichols (braided) Lie algebras of  diagonal type} 
\author{Weicai Wu, Jing Wang, Shouchuan Zhang and  Yao-Zhong Zhang}                 
\lastname{Wu, Wang, Zhang and Zhang}  

\msc{16W30,   16G10}    

\keywords{Braided vector space,   Nichols  algebra,   Nichols (braided) Lie algebra, Graph}         

\address{%
Weicai Wu\\               
Department  of Mathematics\\ 
Zhejiang University\\
Hangzhou 310007\\ 
P.R. China \\            
School of Mathematics\\   
Hunan Institute of Science and Technology,
Yueyang 414006\\     
P.R.  China \\
weicaiwu@hnu.edu.cn                
}
\address{%
Jing Wang\\          
College of Science\\
Beijing Forestry University \\
Beijing 100083\\ 
P.R. China\\
}        
\address{%
Shouchuan Zhang\\               
Department  of Mathematics\\ 
Hunan University\\
Changsha  410082\\ 
P.R. China \\            
sczhang@hnu.edu.cn               
}
\address{%
Yao-Zhong Zhang\\               
School of Mathematics and Physics\\ 
The University of Queensland\\
Brisbane 4072\\ 
Australia \\            
yzz@maths.uq.edu.au               
}

\begin{document}

\maketitle 

\begin {abstract} Let  $V$ be a braided vector space of diagonal type.
Let $\mathfrak B(V)$, $\mathfrak L^-(V)$ and $\mathfrak L(V)$ be the Nichols algebra, Nichols Lie algebra
and Nichols braided Lie algebra over $V$, respectively. We show that a monomial belongs to $\mathfrak L(V)$ if and only if
this monomial is connected. We obtain the basis for $\mathfrak L(V)$ of arithmetic root systems and the dimension of $\mathfrak L(V)$ of finite Cartan type.
We give the sufficient and necessary conditions for  $\mathfrak B(V) = F\oplus \mathfrak L^-(V)$ and $\mathfrak L^-(V)= \mathfrak L(V)$.
We obtain an explicit basis for $\mathfrak L^ - (V)$ over the quantum linear space $V$ with $\dim V=2$.
\end {abstract}



\section {Introduction}\label {s0}
Nichols algebras have found significant applications in various areas of mathematics and mathematical physics including the
theories of pointed Hopf algebras and logarithmic quantum fields. In \cite {He05}, one-to-one correspondences
between Nichols algebras of diagonal type and arithmetic root systems as well as between generalized Dynkin diagrams
and twisted equivalence classes of arithmetic root systems were established.
The problem of finite-dimensionality of Nichols algebras forms a substantial part of the recent investigations
(see e.g. \cite {AHS08,  AS10,  He05,  He06a,  He06b,  WZZ15a,  WZZ15b}. Lie algebra arising from a Nichols algebra was studied in \cite {AAB16}.

Braided Lie algebras  were studied in \cite {Ka77, BMZP92, Gu86, GRR95, Kh99, Pa98, Sc79, BFM96, ZZ03}. The current paper will focus on
Nichols  Lie and braided Lie algebras.  In \cite {He05} and \cite {He06a}, a classification on braided vector spaces of diagonal type
with finite-dimensional Nichols algebras was given. In \cite {WZZ15b, WZZ16}, we studied the relationship between Nichols algebras and Nichols braided Lie algebras.
It was proven that a Nichols algebra is finite-dimensional if and only if the corresponding Nichols (braided) Lie algebra is finite-dimensional.
This provides a new method for determining when a Nichols algebra is finite dimensional.

In this work we will show that a monomial belongs to Nichols braided Lie algebra $\mathfrak L(V)$ of braided vector space $V$ of diagonal type
if and only if this monomial is connected. This is one of the main results in this paper, which enables us to obtain the bases for the Nichols braided Lie algebras of arithmetic root systems and the dimensions of the Nichols braided Lie algebras of finite Cartan type.
We give the sufficient and necessary conditions for  $\mathfrak B(V) = F\oplus \mathfrak L^-(V)$ and $\mathfrak L^-(V)= \mathfrak L(V)$, where
$\mathfrak B(V)$, $\mathfrak L^-(V)$ and $\mathfrak L(V)$ denote Nichols algebra, Nichols Lie algebra
and Nichols braided Lie algebra over $V$, respectively. We also obtain an explicit basis for $\mathfrak L^ - (V)$ over the quantum linear space $V$ with $\dim V=2$.

This paper is organized as follows. In the remaining part of this section we provide some preliminaries and set our notations.
In  Section 2 we show that a monomial belongs to $\mathfrak L(V)$ if and only if this monomial is connected when $V$ is  a  braided vector space of diagonal type.
Section 3 presents a basis for $\mathfrak L(V)$, section 4 gives the basis for $\mathfrak L(V)$ of
arithmetic root systems and  obtains  the dimension of $\mathfrak L(V)$ of finite Cartan type,
and section 5 provides some non-zero monomials for the Nichols algebras $\mathfrak B(V)$. In section 6 we present the sufficient and necessary conditions for
$\mathfrak B(V) = F\oplus \mathfrak L^-(V)$ and $\mathfrak L^-(V)= \mathfrak L(V)$, respectively.
We also give an explicit basis for  $\mathfrak L^ - (V)$ over the quantum linear space $V$ with $\dim V=2$.

\section*{ Preliminaries}

For any matrix $(q_{ij})_{n \times n}$  over $F^*$, define a bicharacter $\chi$ from $\mathbb Z^n \otimes \mathbb Z^n$ to $F^*$ such that $\chi (e_i, e_j) =q_{ij}$  for $1\le i, j \le n$, where $\{ e_1, e_2, \cdots, e_n\}$ is a basis of $\mathbb Z^n.
$ Let $V$ be a vector space  with basis
$x_1,  x_2,  \cdots,  x_n$. Define $\alpha (e_i\otimes  x_j) = q_{ij}x_j$ and $\delta (x_j) =e_j \otimes x_j$ for $1\le i, j \le n$. It is clear that $(V, \alpha, \delta)$ is a Yetter-Drinfeld module over $\mathbb Z^n$ and $(V, C) $ is a  braided vector space under braiding $C,$  where $C(x_i \otimes x_j) = q_{ij}x_j \otimes x_i$
and $C^{-1}(x_i \otimes x_j) = q_{ji} ^{-1}x_j \otimes x_i.$ In this case, $V$ is called a  braided vector space  of diagonal type and $(q_{ij})_{n \times n}$ is called a braiding matrix of $V$.
Throughout this paper  braided vector space $V$ is of diagonal type with basis
$x_1,  x_2,  \cdots,  x_n$ and $C(x_i \otimes x_j) =q_{ij} x_j \otimes x_i$ without special announcement. Let $\mathfrak B(V)$ be the Nichols algebra over the braided  vector space $V$.
Define $p_{ij} := q_{ij}$ for $1\le i, j \le n$  and  $p_{uv} := \chi ( {\rm deg} (u),   {\rm deg }(v))  $ for any homogeneous element $u,   v \in \mathfrak B(V).$
Denote ${\rm ord } (p_{uu})$ the  order of $p_{uu}$ with respect to multiplication. Let $|u|$ denote length of  homogeneous element $u\in \mathfrak B(V).$
Let $D =: \{[u] \mid [u] \hbox { is a hard super-letter }\}$, $\Delta ^+(\mathfrak B(V)): =  \{ \deg (u) \mid [u]\in D\}$,
$ \Delta (\mathfrak B(V)) := \Delta ^+(\mathfrak B(V)) \cup \Delta ^-(\mathfrak B(V))$,
which is called the root system of $V.$
 If $ \Delta (\mathfrak B(V))$ is finite,
then it is called an arithmetic root system. Let  $\mathfrak L(V)$ denote the  braided Lie algebras generated by $V$ in $\mathfrak B(V)$ under Lie operations $[x, y]=yx-p_{yx}xy$,  for any homogeneous elements $x, y\in \mathfrak B(V)$. $(\mathfrak L(V), [\ ])$ is called Nichols braided Lie algebra of $V$. Let  $\mathfrak L^-(V)$ denote the  Lie algebras generated by $V$ in $\mathfrak B(V)$ under Lie operations $[x, y]^-=yx-xy$,  for any homogeneous elements $x, y\in \mathfrak B(V)$. $(\mathfrak L^-(V), [\ ]^-)$ is called Nichols  Lie algebra of $V$.
The other notations are the same as in \cite {WZZ15a}.

\vskip.1in
Recall the dual $\mathfrak B(V^*) $ of Nichols algebra $\mathfrak B(V) $ of rank $n$
in \cite [Section 1.3]{He05}. Let $y_{i}$ be a dual basis of $x_{i}$. $\delta (y_i)=g_i^{-1} \otimes y_i$,  $g_i\cdot y_j=p_{ij}^{-1}y_j$ and $\Delta (y_i)=g_i^{-1}\otimes y_i+y_i\otimes 1.$ There exists a bilinear map
$<\cdot, \cdot >$: $(\mathfrak B(V^*)\# FG)\times\mathfrak B(V)$ $\longrightarrow$ $\mathfrak B(V)$ such that

$<y_i, uv>=<y_i, u>v +g_i^{-1}.u<y_i, v>$ and $<y_i, <y_j, u>>=<y_iy_j, u>$

\noindent for any $u, v\in\mathfrak B(V)$. Furthermore,  for any $u\in \oplus _{i=1}^\infty \mathfrak B(V)_{(i)}$,
one has  $u=0$ if and only if $<y_i, u>=0$ for any $1\leq i \leq n.$
We have
\begin {eqnarray}\label {e0.1} [[u, v], w]=[u, [v, w]]+p_{vw}^{ -1} [[u, w], v]+(p_{wv}-p_{vw}^{-1})v\cdot[u, w],
\end {eqnarray}
\begin {eqnarray}\label {e0.2} [u, v \cdot w]=p_{wu}[uv]\cdot w + v \cdot [uw].\end {eqnarray}

We now recall some basic concepts of the graph theory (see \cite {Ha69}).
Let $\Gamma _1$ be a non-empty set and $\Gamma _2 \subseteq \{  \{ u, v\} \mid u, v \in \Gamma _1, \hbox { with } u \not= v \} \subseteq 2 ^{\Gamma_1}.$ Then $\Gamma = (\Gamma _1, \Gamma_2)$ is called a graph; $\Gamma_1$ is called the vertex set of  $\Gamma$; $\Gamma_2$ is called the edge set of  $\Gamma$; Element $\{u, v\} \in \Gamma_2$ is called an edge, written $a_{u, v}$.
 If $G = (G_1, G_2)$ is a graph and $G_1 \subseteq \Gamma _1$ and $G_2 \subseteq \Gamma _2$, then $G$ is called a subgraph of $\Gamma.$ If $ \emptyset \not= H_1 \subseteq \Gamma_1$ and $H_2 = \{ a_{u, v} \in \Gamma _2\mid u, v \in H_1 \}$, then $H=(H_1, H_2)$ is a subgraph, called the subgraph generated by $H_1$ in $\Gamma$.

 $a_{u_mu_{m-1}} \cdots a_{u_3u_2} a_{u_2u_1}$ is called a walk  from $u_1$ to $u_m.$
 We can define an equivalent relation on $\Gamma_1$ as follows: for any $u, v\in \Gamma _1$, $u$ and $v$ are equivalent if and only if there exists a walk from $u$ to $v$ or $u=v.$  Every subgraph  generated by  every equivalent class of $\Gamma_1$ is called
a connected component of $\Gamma.$

Let $\Gamma (V)$ be the  generalized Dynkin diagram of $V$ with $p_{x_i,  x_i}$  and $\widetilde{p}_{x_i,  x_j}$ for $1\le i\not= j \le n $ omitted. This is called a pure generalized Dynkin graph of $V$, i.e. $\Gamma (V)_1 = \{ 1, 2, \cdots, n\}$ and $\Gamma (V)_2 = \{  a_{ij} \mid  p_{ij} p_{ji} \not= 1, i \not= j\}$.

Let   $u=  h_1 h_2\cdots  h_m$ be   a  monomial with $ h_j = x_{i_j}$  for $1\leq i_1,i_2,\ldots,i_m\leq n,1\le j \le m$. ${\rm deg} (u) = \lambda _1 e_1 + \cdots + \lambda _n e_n, $ where $ {\rm deg }(x_i) = e_i.$ Let ${\rm deg}_{x_i} (u) := \lambda _i$. Let $\mu (u):= \{ x_{i_1},  \cdots,  x_{i_m}\}$  and  $\Gamma (u)$ be a pure generalized Dynkin subgraph generated by $\mu (u)$.
 If $\Gamma (u)$
is connected,  then $u$ is called connected (or $\mu (u) $ is called connected). Otherwise, $u$ is called disconnected (or $\mu (u) $ is called disconnected).

For $u, v \in \mathfrak B(V)$, if there exists a non-zero $a\in F$ such that  $u = av$, then we write  $u \sim v.$ This is an equivalent relation.
If there exist $x_i \in \mu (u)$ and   $x_j \in \mu (v)$   such that $\widetilde{p} _{x_i,  x_j} \not=1$,  then we say that it is connected between monomial $u$ and monomial $v$, written $u \diamondsuit v$ in short. Otherwise, we say that it is disconnected between monomial $u$ and monomial $v$.

Remark:
When $u\not=0$, $\mu (u)$ is independent of the choice of $h_1, h_2, \cdots, h_m,$ since $ \mathfrak B(V)$ is graded and ${\rm deg } (u) $ is unique. Therefore, the connectivity of $u$ is independent of the choice of $h_1, h_2, \cdots, h_m.$ When $u=0$, $\mu (u)$ is dependent on the choice of $h_1, h_2, \cdots, h_m.$  For example, ${\rm deg } (V) >1$, $p_{12} p_{21} =1, $ $p_{11} =-1,$ $u= x_1 ^2 = x_1^2 x_2=0.$ Therefore,  $\mu (u) =\{ x_1\}$  and $u$ is connected. Meantime, $\mu (u) =\{ x_1, x_2\}$  and $u$ is disconnected.

Throughout,  $\mathbb Z =: \{x \mid  x \hbox { is an integer}\}.$ $\mathbb R =: \{ x \mid x \hbox { is a real number}\}$.
$\mathbb N_0 =: \{x \mid  x \in \mathbb Z, x\ge 0\}.$
$\mathbb N =: \{x \mid  x \in \mathbb Z,  x>0\}$.  $F$ denotes the base field,   which is an algebraic closed field with  characteristic zero. $F^{*}=F\backslash\{0\}$. $\mathbb S_{n}$ denotes symmetric group, $n\in\mathbb N$. For any set $X$, $\mid X\mid$ is the cardinal of $X$. ${ {\rm int} }(a)$ means the biggest integer not greater than $a\in \mathbb R$.

\section{Structure of Nichols braided Lie algebras}\label {s1}

In this section we prove that a monomial belongs to $\mathfrak L(V)$ if and only if
this monomial is connected when $V$ is  a  braided vector space of diagonal type.

\begin{Lemma}\label{2.1} Assume that $u, v, w$ are homogeneous elements in $\mathfrak L(V)$. Let $a,  b,  c,  d,  e$ and $f$ denote  $1 - p_{wv}p_{vw},  1 -p_{uw}p_{wu},  1 - p_{uv}p_{vu}$,  $1 - p_{uv}p_{vu}p_{uw}p_{wu}$,  $1 - p_{uv}p_{vu}p_{wv}p_{vw}$ and  $1 - p_{wv}p_{vw}p_{uw}p_{wu}$, respectively. If $\mid\{r|r\in\{uv,uw,vw\} \hbox { and }r\in\mathfrak L(V)\}\mid \geq 2$ and  $ \mid \{t|t\in\{a, b, c\} \hbox { and } t\neq0\} \mid \ge  1$,  then
 $uvw, uwv, vwu, vuw, wuv$, $wvu$  $\in\mathfrak{L}(V)$.
\end {Lemma}
\noindent {\it Proof.} Without loss of generality, we let $uv,uw\in\mathfrak L(V)$. If $vuw\notin\mathfrak L(V)$, then $e=0,f=0$ by
\cite [Lemma 3.2(i)] {WZZ15b} and \cite [Lemma 4.12] {WZZ15a}. If $a=0$, then $b=0,c=0,$  which  is a contradiction. If $a\neq0$, then $b\neq0,c\neq0$, one obtains a contradiction to \cite [Lemma 3.1] {WZZ15b}.  Consequently,   $vuw\in\mathfrak L(V)$. Similarly,  we can obtain others.  \hfill $\Box$

\begin {Lemma} \label {2.2} Assume that $h_i \in \{x_1,  \cdots,  x_n\}$ for $1\le i \le m.$

{\rm (i)} If it is disconnected between monomial $u$ and monomial $v$ (i.e. $\widetilde{p} _{x_i,  x_j} =1$  for any  $x_i \in \mu (u)$,  $x_j \in \mu (v)$ ),  then $[u,  v] =0.$

{\rm (ii)} If $\mu(h_{1}h_{2}\cdots h_{m})$ is  disconnected,  then $\sigma (h_{1},  h_{2},  \cdots,   h_{m})=0$ for any  method $\sigma$ of adding bracket on $h_1,  h_2,  \cdots,  h_m.$

{\rm (iii)} If $h_{1}h_{2}\cdots h_{m}\neq0$ and $\mu(h_{1}h_{2}\cdots h_{m})$ is disconnected,  then $h_{1}h_{2}\cdots h_{m}\notin\mathfrak{L}(V)$.

\end {Lemma}

\noindent {\it Proof.} {\rm (i)}  $u$ and $v$ are quantum commutative
(i.e. $uv= p_{u, v}vu$ ) since $x_i x_j = p_{x_i,  x_j}x_jx_i$ for any  $x_i \in \mu (u)$,  $x_j \in \mu (v)$.

{\rm (ii)}
We show this by induction on $m$. $ [h_1,  h_2] =0$ for $m=2.$  For $m>2,  $ $\sigma (h_{1},  h_{2},  \cdots,  h_{m}) $

\noindent $= [\sigma _1 ( h_1 h_2 \cdots h_t ),  \sigma _2 ( h_{t+1} h_{t+2} \cdots h_m ]$.
If both $ h _1h_2\cdots h_t $ and $ h_{t+1} h_{t+2} \cdots h_m $ are connected,  then it is disconnected between  $h_1 h_2 \cdots h_t $ and $ h_{t+1} h_{t+2} \cdots h_m $. By Part {\rm (i)},   $\sigma (h_{1},  h_{2},  \cdots,  h_{m})=0$.
 If either  $h_1h_2 \cdots h_t $ or $ h_{t+1} h_{t+2} \cdots h_m$ is  disconnected,
then either  $\sigma _1 (h_1h_2\cdots h_t )=0$ or $ \sigma _2 (h_{t+1} h_{t+2} \cdots h_m)=0 $ by inductive hypothesis.

{\rm (iii)}  It follows from Part {\rm (ii)}.
\hfill $\Box$

\begin {Lemma} \label {2.3}

{\rm (i)} If $u$ is  a monomial, then there exist monomials $u_{1},u_{2},\ldots,u_{r}$ such that $u \sim  u_1u_2 \cdots u_r $ and  $\{ \Gamma (u_1),  \Gamma ( u_2),  \cdots,  \Gamma(u_r) \}$ is complete set of  connected components of $\Gamma (u)$ with ${\rm deg}( u) =\sum \limits _{i=1} ^r {\rm deg}( u_i)$   (which is called a  decomposition of  connected components  of $u$).

{\rm (ii)}  If a monomial $u$ is connected with $\mid u \mid >1$,  then there exist two  connected  monomials $v$ and $w$ such that $v\diamondsuit  w$  with $u \sim vw$.

\end {Lemma}
\noindent {\it Proof.} {\rm (i)}  If $u$ is connected, it is clear. If $u$ is not connected, we show it by induction on $\mid u \mid $. It is clear when  $\mid u \mid =1.$  Now assume $\mid u \mid >1$. Let $\Omega := \{ v \mid
    v \hbox { is a connected monomial  }  \}$.  Let $ u \sim  v_1v_2v_3$ and $v_2$ be in $\Omega$ such that $\mid v_2 \mid = max \{ \mid v \mid \ \  \mid  v \in  \Omega, u\sim v_1v v_3,   v_1 \hbox { and } v_3 \hbox { are monomials or  } 1 \}$.

    Now we show that $v_2$ is a connected component  of $u$. In fact, if it does not hold, then $v_1 \diamondsuit v_2$ or $v_3 \diamondsuit v_2$. Without the lose of generality ¡¡
We assume $v_1 \diamondsuit v_2$. Then there exists $t$ such that $ h_{i_t}\diamondsuit  v_2$ with
$v_1 = h _{i_1} \cdots  h_{i_s} $  and  $v_2  h_{i_j} \sim  h_{i_j } v_2$ for $t+1< j \leq s$. Consequently,  $u \sim w_1 h_{i_t}v_2 w_3 $ and $h_{i_t}v_2 $ is connected. which is a contradiction.

We have obtained that it is disconnected between $v_1$ and $  v_2$ in proof above.
Thus  $u =v_1 v_2 v_3 \sim v_2 v _1v_3 $ and $v_1v_3 \sim u_2u_3 \cdots u_r$ is a  decomposition of  connected components  of $v_1v_3$ since $\mid v_1v_3 \mid < \mid u \mid. $ Setting  $u_1:= v_2 $ we complete  proof of Part {\rm (i)}.

{\rm (ii)} Let $ u = h_1 h_2\cdots  h_m$   and $ h_2\cdots  h_m \sim u_1u_2\cdots u_r$ be a decomposition of connected components of $ h_2\cdots  h_m$.
Consequently,  $ h_1u_1\cdots u_{r-1}$ is connected since $ h_1\diamondsuit  u_i$ for $1\le i \le r$.
\hfill  $\Box$

\begin {Theorem} \label {2.4} Assume that $h_i \in \{x_1,  \cdots,  x_n\}$ and $p_{h_1,   h_1} \not= 1$ when   $\mu (h_{1}\cdots  h_{m})$ $= h_1$.

{\rm (i)}
 If $ h_{1}\cdots  h_{b},  h_{b+1}\cdots  h_{m}\in\mathfrak{L}(V)$ and
$ h_{1}\cdots  h_{b}\diamondsuit  h_{b+1}\cdots  h_{m}$ with $h_{1}\cdots  h_{b}\not= 0$ or $ h_{b+1}\cdots  h_{m} \not= 0$,
then there exist $\tau\in \mathbb  S_m $ such that $ h_{1}\cdots  h_{m}\sim h_{\tau(1)}\cdots  h_{\tau(m)}$ with

\noindent $ h_{\tau(1)}\cdots h_{\tau(m-1)}\in\mathfrak{L}(V)$ and  $ h_{\tau(1)}\cdots h_{\tau(m-1)}\diamondsuit  h_{\tau(m)}$,
or with $ h_{\tau(2)}\cdots  h_{\tau(m)}\in\mathfrak{L}(V)$ and $ h_{\tau(2)}\cdots  h_{\tau(m)}\diamondsuit  h_{\tau(1)}$.

{\rm (ii)}
If $  h_{1} h_{2}\cdots  h_{m-1}\in\mathfrak{L}(V) $  and $  h_{1} h_{2}\cdots  h_{m-1}$ $\diamondsuit  h_m$ or  $   h_{2} h_{3}\cdots h_{m}\in\mathfrak{L}(V) $  and $  h_{2} h_{3}\cdots  h_{m}$ $\diamondsuit  h_1, $ then $ h_{1} h_{2}\cdots  h_{m}\in\mathfrak{L}(V)$.

{\rm (iii)} If $0\not=   h_{1} h_{2}\cdots  h_{m}\in\mathfrak{L}(V) $,  then there exists $\tau \in \mathbb S_m$ such that $ h_{1}\cdots  h_{m}\sim h_{\tau(1)}\cdots  h_{\tau(m)}$ with $ 0\neq h_{\tau (1)} h_{\tau(2)}\cdots h_{\tau(m-1)}\in\mathfrak{L}(V)$ and $ h_{\tau (1)}$ $h_{\tau(2)}$ $\cdots$  $h_{\tau(m-1)}$ $\diamondsuit$ $h_{\tau(m)}$,   or $0\neq h_{\tau (2)} h_{\tau(3)}\cdots  h_{\tau(m)}\in\mathfrak{L}(V)$ and $ h_{\tau (2)}$ $h_{\tau(3)}$ $\cdots$ 
$h_{\tau(m)}$ $\diamondsuit$ $ h_{\tau(1)}$ .

 {\rm (iv)}
If monomial $u= h_1h_2 \cdots h_m$ is connected,  then $u \in  \mathfrak{L}(V). $

\end {Theorem}
\noindent {\it Proof.}
We show {\rm (i)},  {\rm (ii)},  {\rm (iii)} and {\rm (iv)} by induction on the length of $| h_{1} h_{2}\cdots  h_{m}|=m$.

  Assume $m=2$. {\rm (i)} and {\rm (iv)} are clear. If $ h_{1}\neq  h_{2}$,  then {\rm (ii)}, {\rm (iii)} follows from \cite [Lemma 4.12]{WZZ15a}, \cite [Lemma 5.2]{WZZ15a},  respectively. If $ h_{1}=  h_{2}$,  then {\rm (ii)} and {\rm (iii)} follow from \cite [Lemma 4.3]{WZZ15a}  and  \cite [Lemma 1.3.3{\rm (i)}]{He05}.

Now  $m>2$.

{\rm (i)} If $b=m-1$ or $b=1$,  let $\tau=id$. Now assume that $1<b<m-1$. There exist $1 \le a \le b$ and $b+1 \le c \le m$ such that $\widetilde{p} _{h_a,  h_c} \not= 1$.

{\rm (1)}. Assume $ h_{1}\cdots  h_{b}\neq0$. We show this by  induction on $b$.
There exist $\tau\in \mathbb  S_{\{1, \ldots, b\}}$ such that $ h_{1}\cdots  h_{b}\sim  h_{\tau(1)}\cdots  h_{\tau(b)}$ with
 case (a):  $0\neq h_{\tau(1)}\cdots h_{\tau(b-1)}\in\mathfrak{L}(V)$, $ h_{\tau (1)} h_{\tau(2)}\cdots  h_{\tau(b-1)}\diamondsuit  h_{\tau(b)}$
or with case (b): $0\neq  h_{\tau(2)}\cdots h_{\tau(b)}\in\mathfrak{L}(V)$,  $ h_{\tau (2)} h_{\tau(3)}\cdots  h_{\tau(b)}\diamondsuit  h_{\tau(1)}$ by induction hypotheses of {\rm (iii)}.

(a)$_{1}$.  If $ h_{b+1} \cdots  h_{m}\diamondsuit  h_{\tau(b)}$,  then
$ h_{\tau(b)} h_{b+1}\cdots  h_{m}\in\mathfrak L( V)$ by induction hypotheses of {\rm (ii)},  and $h_{\tau(1)}\cdots h_{\tau(b-1)}\diamondsuit  h_{\tau(b)} h_{b+1}\cdots  h_{m}$ since
$ h_{\tau (1)} h_{\tau(2)}\cdots  h_{\tau(b-1)}\diamondsuit  h_{\tau(b)}$. It is proved since
$| h_{\tau(1)}\cdots  h_{\tau(b-1)}|=b-1<b$ and induction hypotheses.

(a)$_{2}$.  If  $\tilde{p}_{ h_{\tau(b)},  h_{j_{1}}}=1$ for all $j_{1}\in\{b+1, \ldots, m\}$,  then
$ h_{\tau(b)} h_{b+1}\cdots  h_{m}\sim  h_{b+1}\cdots  h_{m} h_{\tau(b)}$. On the other hand,  $h_{\tau(1)}\cdots h_{\tau(b-1)}\diamondsuit   h_{b+1}\cdots  h_{m}$ since
$ h_{1}\cdots  h_{b}\diamondsuit $ $ h_{b+1}\cdots  h_{m}$.
 Consequently, $ h_{\tau_{1}\tau(1)}\cdots  h_{\tau_{1}\tau(b-1)} h_{\tau_{1}(b+1)}\cdots h_{\tau_{1}(m)}\in\mathfrak L( V)$ and   $ h_{\tau(1)}$ $\cdots$  $h_{\tau(b-1)}$ $h_{b+1}$ $\cdots$  $h_{m} $ $\sim  h_{\tau_{1}\tau(1)}$ $\cdots $ $h_{\tau_{1}\tau(b-1)} h_{\tau_{1}(b+1)}\cdots  h_{\tau_{1}(m)}$ for some $\tau_{1}\in \mathbb  S_{\{\tau(1), \ldots, \tau(b-1), b+1, \ldots, m\}} $ by the induction hypotheses of {\rm (i)} and {\rm (ii)}. 
 $h_{\tau_{1}\tau(1)}$ $\cdots $ $h_{\tau_{1}\tau(b-1)}$ $ h_{\tau_{1}(b+1)}$ $\cdots$ $ h_{\tau_{1}(m)}$ $\diamondsuit$ $  h_{\tau(b)}$  since $ h_{\tau (1)} h_{\tau(2)}\cdots  h_{\tau(b-1)}\diamondsuit  h_{\tau(b)}$. See
 \begin {eqnarray*}
 h_{1}\cdots  h_{m}
&\sim & h_{\tau(1)}\cdots  h_{\tau(b)} h_{b+1}\cdots  h_{m}
\\
&\sim & h_{\tau(1)}\cdots h_{\tau(b-1)} h_{b+1}\cdots  h_{m} h_{\tau(b)} \\
&\sim & h_{\tau_{1}\tau(1)}\cdots  h_{\tau_{1}\tau(b-1)} h_{\tau_{1}(b+1)}\cdots h_{\tau_{1}(m)} h_{\tau(b)}\\
&=&  h_{\tau_{1}\tau(1)} \cdots h_{\tau_{1}\tau(b-1)}h_{\tau_{1}\tau(b+1)}\cdots  h_{\tau_{1}\tau(m)} h_{\tau_{1}\tau(b)}
 \end {eqnarray*}
 by $ h_{\tau_{1}(b+1)}\cdots  h_{\tau_{1}(m)}$
$= h_{\tau_{1}\tau(b+1)}\cdots  h_{\tau_{1}\tau(m)}$ since $\tau\in \mathbb  S_{ \{1, 2, \ldots, b\}}$,
 and
$ h_{\tau_{1}\tau(b)}= h_{\tau(b)}$ since $\tau_{1}\in \mathbb  S_{\{\tau(1), \tau(2), \ldots, \tau(b-1), b+1, \ldots, m \}}$. Consequently,
 $$\tau':={\tiny \left(       \begin{array}{ccccccc}
              1&\cdots&b-1&b&\cdots&m-1&m \\
               \tau_{1}\tau(1)&\cdots&\tau_{1}\tau(b-1)&\tau_{1}\tau(b+1)&\cdots& \tau_{1}\tau(m)&\tau_{1}\tau(b)
                   \end{array}\right) } \in \mathbb  S_{m}, $$
                   $ h_{\tau'(1)}\cdots h_{\tau'(m-1)}\in\mathfrak{L}(V), h_{\tau'(1)}\cdots h_{\tau'(m-1)}\diamondsuit  h_{\tau'(m)}$ and
$ h_{1}\cdots  h_{m}\sim h_{\tau'(1)}\cdots  h_{\tau'(m)}$.

(b)$_{1}$.  If $ h_{\tau(2)}\cdots h_{\tau(b)}\diamondsuit  h_{b+1}\cdots  h_{m}$,  then  $ h_{\tau_{1}\tau(2)}\cdots  h_{\tau_{1}\tau(b)} h_{\tau_{1}(b+1)}\cdots h_{\tau_{1}(m)}\in\mathfrak L( V)$  and $ h_{\tau(2)}\cdots h_{\tau(b)} h_{b+1}\cdots  h_{m}
\sim  h_{\tau_{1}\tau(2)}\cdots
 h_{\tau_{1}\tau(b)} $ \ $h_{\tau_{1}(b+1)}\cdots  h_{\tau_{1}(m)}$ for some $\tau_{1}\in \mathbb  S_{\{\tau(2), \ldots, \tau(b), b+1, \ldots, m\}}$ by the induction hypotheses of {\rm (i)} and {\rm (ii)}. We know $h_{\tau_{1}\tau(2)}\cdots  h_{\tau_{1}\tau(b)} h_{\tau_{1}(b+1)}\cdots h_{\tau_{1}(m)}\diamondsuit  h_{\tau(1)}$ since $ h_{\tau (2)} h_{\tau(3)}\cdots  h_{\tau(b)}\diamondsuit  h_{\tau(1)}$.
We obtain
$ h_{1}\cdots  h_{m}\sim  h_{\tau(1)}\cdots  h_{\tau(b)} h_{b+1}\cdots h_{m}\sim  h_{\tau(1)} h_{\tau_{1}\tau(2)} $ $\cdots h_{\tau_{1}\tau(b)} h_{\tau_{1}(b+1)} $ $\cdots  h_{\tau_{1}(m)}= h_{\tau_{1}\tau(1)} h_{\tau_{1}\tau(2)}$ $\cdots  h_{\tau_{1}\tau(b)} $ $h_{\tau_{1}\tau(b+1)}\cdots h_{\tau_{1}\tau(m)}$ by

\noindent $ h_{\tau_{1}(b+1)}\cdots h_{\tau_{1}(m)}= h_{\tau_{1}\tau(b+1)}\cdots  h_{\tau_{1}\tau(m)}
$ since
$\tau\in \mathbb  S_{ \{1, \ldots, b\}}$ and $ h_{\tau_{1}\tau(1)}= h_{\tau(1)}$ since
$\tau_{1}\in \mathbb  S_{\{\tau(2), \ldots, \tau(b), b+1, \ldots, m\}}$. Therefore,  $\tau':=\tau_{1}\tau \in \mathbb S_m$,   $ h_{\tau'(2)}\cdots h_{\tau'(m)}\in\mathfrak{L}(V),$ $ h_{\tau'(2)}\cdots h_{\tau'(m)}\diamondsuit  h_{\tau'(1)}$ and $ h_{1}\cdots  h_{m}\sim  h_{\tau'(1)}\cdots  h_{\tau'(m)}$.

(b)$_{2}$. If $\tilde{p}_{ h_{\tau(j_{1})},  h_{j_{2}}}=1$ for all $j_{1}\in\{2, \ldots, b\}, j_{2}\in\{b+1, \ldots, m\}$,  then
$ h_{\tau(2)}\cdots  h_{\tau(b)} $ $h_{b+1}\cdots  h_{m}\sim  h_{b+1}\cdots h_{m} h_{\tau(2)}\cdots  h_{\tau(b)}
$. On the other hand, $h_{\tau(1)}$ $\diamondsuit$ $  h_{b+1}$ $\cdots$ $  h_{m}$ since $ h_{1}\cdots  h_{b}$ $\diamondsuit$  $h_{b+1}\cdots  h_{m}$.   Then
$ h_{\tau(1)} h_{b+1}\cdots  h_{m}\in\mathfrak L( V)$ by induction hypotheses of {\rm (ii)}. See $ h_{1}$ $\cdots$ 
$  h_{m}$ $\sim $ $ h_{\tau(1)}$ $\cdots$ $  h_{\tau(b)}$ $ h_{b+1}$ $\cdots$ $  h_{m}\sim  h_{\tau(1)} h_{b+1}\cdots  h_{m} h_{\tau(2)}\cdots h_{\tau(b)}$, and $h_{\tau(1)} h_{b+1}\cdots  h_{m} \diamondsuit  h_{\tau(2)}\cdots h_{\tau(b)}$ since $h_{\tau(1)}$ $\diamondsuit$ $  h_{\tau(2)}$ $\cdots$ $ h_{\tau(b)}$.
{\rm (i)} holds since
$| h_{\tau(2)}\cdots  h_{\tau(b)}|=b-1<b$ and induction hypotheses on $b$.

{\rm (2)}. If  $ h_{1}\cdots  h_{b}= 0$,  then  $ h_{b+1}\cdots  h_{m}\neq0$ and $ h_{b+1}\cdots  h_{m}$ is connected by Lemma \ref {2.2}{\rm (iii)}. Consequently,   $0= h_{1}\cdots  h_{m-1}  \in \mathfrak L(V)$ and $h_{1}\cdots  h_{m-1}\diamondsuit  h_m.$

{\rm (ii)} Assume that $ h_{1} h_{2}\cdots  h_{m}\notin\mathfrak L( V)$. Obviously,
$h_{1} h_{2}\cdots  h_{m-1} \not= 0$ and $h_{2} h_{3}\cdots  h_{m} \not=0.$
Set $i _0 :=1$ and   $j_0 := m-1$ when $h_{1}h_{2}\cdots h_{m-1}\in\mathfrak L( V)$;
 $i_0 :=2$ and   $j _0 := m$ when $h_{2}h_{3}\cdots h_{m}\in\mathfrak L( V)$;
$N_0 := \{ i_0, i_0+1, \cdots,  j_0\};  $ $A_0   :=\{1, \ldots, m\}-N_0$,

Now we prove the following $\mathbf{Assertion (k )}$ by induction on $k$,  $1\le k \le m-2.$

\noindent
$\mathbf{Assertion (k ):}$ There exist $1\le i_k \le  j_k\le m$,  $\tau _k  \in \mathbb S_{N_{k-1}}$ such that the following conditions hold:
${\rm   (C_1)}.$ $0\neq h_{\tau^{k}(i_k)}\cdots h_{\tau^{k}(j_k)}\in\mathfrak L( V);$
${\rm   (C_2)}.$ $\tilde{p}_{h_{r}, h_{t}}=1$ for $\forall \ r\neq t\in A_k$; ${\rm   (C_3)}.$
$h_{\tau^{k}(i_k)}\cdots h_{\tau^{k}(j_k)}\diamondsuit  h_r$ for  $\forall\ r\in A_k$; $ {\rm   (C_4)}.$
$\tilde{p}_{ h_{\tau^{k}(i_k)}\cdots h_{\tau^{k}(j_k)}   ,  h_{r}}=1$  for  $\forall\ r\in A_k$; ${\rm   (C_5)}.$ $h_{\tau^{k}(i_{k-1})}h_{\tau^{k}(i_{k-1}+1)}\cdots h_{\tau^{k}(j_{k-1})} \sim  h_{\tau^{k-1}(i_{k-1})}h_{\tau^{k-1}(i_{k-1}+1)}\cdots h_{\tau^{k-1}(j_{k-1})}  $
 and  $h_{\tau^{k}(1)} \cdots  h_{\tau^{k}(m)} \sim   h_1 h_2 \cdots h_m$; ${\rm   (C_6)}.$
 $i_{k-1} \le i_k\le j_{k} \le j_{k-1}$ with  $j_{k-1}- i_{k-1} = j_k -i_k +1$,  where  $N_k := \{\tau^{k} ({i_k}),  \tau^{k}({i_k+1}),  \cdots,  \tau^{k}({j_k})\} $
and  $A_k   :=\{1, \ldots, m\}-N_k$,  $\tau^{k}:=\tau_{k}\tau_{k-1}\cdots \tau_{1}$.

 Step 1. For $k=1, $
 now we construct $i _1$ and $j_1$ as follows.
Considering  $ 0 \not= h_{i_0}\cdots h_{j_0}\in\mathfrak{L}(V)$, we have
 $\tilde{p}_{h_{i_0}\cdots h_{j_0}, h_{r}}=1$  for  $r \in A_0$ by \cite [Lemma 4.12]{WZZ15a} and there exist $\tau_{1}\in \mathbb S_{N_0}$ such that $h_{i_0}\cdots h_{j_0}\sim h_{\tau_{1}(i_0)}\cdots h_{\tau_{1}(j_0)}$ with  ¡¡case (a) $0\neq h_{\tau_{1}(i_0)}\cdots h_{\tau_{1}(j_0-1)}\in\mathfrak{L}(V), $  $ h_{\tau_{1}(i_0)}\cdots h_{\tau_{1}(j_0-1)}  \diamondsuit h_{\tau_{1}(j_0)} $,
  $i_1 := i_0$ and  $j_1 := j_0-1$
or with case  (b) $0\neq h_{\tau_{1}(i_0 +1)}\cdots h_{\tau_{1}(j_0)}\in\mathfrak{L}(V), $
 $ h_{\tau_{1}(i_0 +1)}\cdots h_{\tau_{1}(j_0)}  \diamondsuit h_{\tau_{1}(i_0)} $,
 $i_1 := i_0 +1$ and $j_1 := j_0$ by induction hypotheses of (iii). Obviously,
$0\neq h_{\tau^{1}(i_1)}\cdots h_{\tau^{1}(j_1)}\in\mathfrak{L}(V)$ and
$h_{\tau_{1}(r)}=h_{r}$ for $r\in A_0$.

Obviously,  ${\rm   (C_1)}$,  ${\rm (C_5)}$ and ${\rm (C_6)}$ hold.

$\tilde{p}_{ h_{t},  h_{r}} = 1$ for  any $t \not= r  \in A_1, $ i.e. ${\rm (C_2)}$ holds.  Indeed,  if $\tilde{p}_{ h_{t},  h_{r}}\neq1$ for $t \in A_1- A_0,  $  $r \in A_0$,  then $\tilde{p}_{ h_{i_1}\cdots  h_{j_1},  h_{r}}\neq1$ since $\tilde{p}_{ h_{i_0}\cdots  h_{j_0},  h_{r}}=1$.
We obtain
$ h_{i_1}\cdots  h_{j_1} h_{t} h_{r}, $  $ h_t h_{i_1}\cdots  h_{j_1} h_{r},  $
$h_r h_{i_1}\cdots  h_{j_1} h_{t} $ $
\in\mathfrak L( V)$ by Lemma \ref {2.1}. However,  $h_{1}\cdots  h_{m}$ is
a quantum equivalent with one among
$ h_{i_1}\cdots  h_{j_1} h_{t} h_{r}, $  $ h_t h_{i_1}\cdots  h_{j_1} h_{r},  $
$h_r h_{i_1}\cdots  h_{j_1} h_{t} $, $h_r h_{t}h_{i_1}\cdots  h_{j_1} $, which   contradicts  to
   $h_{1}\cdots  h_{m}\notin\mathfrak L( V)$.

   ${\rm (C_3)}$ and ${\rm (C_4)}$ follow from ${\rm (C_2)}$.

Step 2.  Assuming that $\mathbf{Assertion (k )}$ hold,  we prove that $\mathbf{Assertion (k+1)}$ holds,  $k \leq m-3$.
Let $i:=i_k$ and $j:=j_k$ in this proof for convenience.

Considering $0\neq  h_{\tau^{k}(i)}\cdots  h_{\tau^{k}(j)}\in\mathfrak L( V)$,  we have that  there exists
$\tau_{k+1}\in \mathbb  S_{ N_k} $ such that
$ h_{\tau^{k}(i)}\cdots h_{\tau^{k}(j)}\sim  h_{\tau^{k+1}(i)}\cdots  h_{\tau^{k+1}(j)}$ with case ($1^{\circ }$):
$0\neq $ $ h_{\tau^{k+1}(i)}$ $\cdots$ $  h_{\tau^{k+1}(j-1)}$ $\in$ $\mathfrak L( V)$, $h_{\tau^{k+1}(i)}\cdots  h_{\tau^{k+1}(j-1)}\diamondsuit h_{\tau^{k+1}(j)}$,
  $\tau^{k+1}¡¡ (i_{k+1}) := \tau^{k+1}(i)$  and $\tau^{k+1}¡¡(j_{k+1}):= \tau^{k+1}(j-1)$,
or with case ($2^{\circ }$):
\noindent $0\neq  h_{\tau^{k+1}(i+1)}\cdots  h_{\tau^{k+1}(j)}\in\mathfrak L( V), $
$h_{\tau^{k+1}(i+1)}\cdots  h_{\tau^{k+1}(j)}\diamondsuit h_{\tau^{k+1}(i)}$,
 $ \tau^{k+1}¡¡(i_{k+1} ):= \tau^{k+1}(i+1)$  and $\tau^{k+1}¡¡(j_{k+1}):= \tau^{k+1}(j)$ by induction hypotheses of (iii). We obtain
$ h_{1}\cdots  h_{m}\sim  h_{\tau^{k+1}(1)}\cdots  h_{\tau^{k+1}(m)}$.
For convenience, let  $ h_{s}':= h_{\tau^{k+1}(s)}$ for all $s\in\{1, \ldots, m\}$.
Obviously,  $\tau_{k+1}(t)
=t$ for all $t\in A_k $.

Obviously,  ${\rm (C_1)}$,  ${\rm (C_5)}$ and ${\rm (C_6)}$ hold.

($1^{\circ}$)    $A_{k+1} = A_k \cup \tau ^{k+1} (j)$ and $N_{k+1} = N_k - \tau ^{k+1} (j).$ If there exists $r\in A_{k+1} $ such that $\tilde{p}_{ h_{\alpha}',  h_{r}}= 1$ for all $\alpha\in\{i, \ldots, j-1\}$,  then
$\tilde{p}_{ h_{j}',   h_{r}}\neq 1$ by $\mathbf{Assertion (k )}$ ${\rm (C_3)}$. Therefore  $\tilde{p}_{ h_{i}'\cdots  h_{j}',   h_{r}}\neq 1$,  which contradicts to $\mathbf{Assertion (k )}$ ${\rm (C_4)}$. For all $r\in A_{k+1} $, there exists $\alpha\in\{i, \ldots, j-1\}$ such that
\begin {eqnarray} \label {e2.4.100}\tilde{p}_{ h_{\alpha}',   h_{r}}\neq 1\end {eqnarray} and
 $ h_{1}'\cdots  h_{i-1}' h_{i}'\cdots  h_{j-1}'\in\mathfrak L( V)$ by induction hypotheses of {\rm (iv)} and Lemma \ref {2.2} {\rm (ii)}.

Set $ \{g_1 , \cdots ,  g_{\beta},  l_1 ,  \cdots , l_{\eta}\} = \{\tau^{k+1}(j+1), \ldots, \tau^{k+1}(m-1)\} \subseteq A_k$  with $g_1 < \cdots < g_{\beta},  l_1 < \cdots < l_{\eta}$ such that
$\tilde{p}_{ h_{j}',   h_{g_{\lambda}}}\not= 1$ and $\tilde{p}_{ h_{j}',   h_{l_{\xi}}}= 1$ for all $ \lambda  \in\{1, \ldots, \beta\}$,  all \ $\xi\in\{1, \ldots, \eta\}$.
Then $ h_{1}'\cdots  h_{j-1}' h_{l_{1}}'\cdots  h_{l_{r}}'\in\mathfrak L( V)$
for all $r\in\{1, \ldots, \eta\}$ and $ h_{j}' h_{g_{1}}'\cdots h_{g_{t}}'\in\mathfrak L( V)$
for all $t\in\{1, \ldots, \beta\}$ by induction hypotheses of {\rm (iv)}.
If $\tilde{p}_{ h_{j}',   h_{m}'}\not= 1$,  then
$\tilde{p}_{ h_{1}'\cdots  h_{j-1}' h_{l_{1}}'\cdots  h_{l_{\eta}}',   h_{m}'} $ $\not= 1$ and
$\tilde{p}_{ h_{j}' h_{g_{1}}'\cdots  h_{g_{\beta}}',   h_{m}'}\not= 1$ by $\mathbf{Assertion (k )}$ ${\rm C_2 }$ and ${ \rm  C_4 }$. Then $ h_{1}'\cdots  h_{j-1}' h_{l_{1}}'\cdots h_{l_{\eta}}' h_{j}' h_{g_{1}}'\cdots  $ $ h_{g_{\beta}}' h_{m}'
\in  \mathfrak L( V)$ by Lemma \ref {2.1},
$ h_{1}\cdots  h_{m}\sim  h_{1}'\cdots  h_{m}' \sim   h_{1}'\cdots h_{j-1}' h_{l_{1}}'\cdots  h_{l_{\eta}}' h_{j}' h_{g_{1}}'\cdots  h_{g_{\beta}}' h_{m}'$. It is a contradiction. Thus  $\tilde{p}_{ h_{j}',   h_{m}'}= 1$ and                                                                                                                                                                                                                                                                                                                                                                                                                                                                                                                                                                                                                                                                                                                                                                                                                                                                                                                                                                                                                                                                                                                                                                                                                                                                                                                                                                                                                                                                                                                                                                                                                                                                                                                                                                                                                                                                                                                                                                                                                                                                                                                                                                                                                                                                                                                                                                                                                                                                                                                                                                                                                                                                                                                                                                                                                                                                                                                                                                                                                                                                                                                                                                                                                                                                                                                                                                                                                                                                                                                                                                                                                                                                                                                                                                                                                                                                                                                                                                                                                                                                                                                                                                                                                                                                                                                                                                                                                                                                                                                                                                                                                                                                                                                                                                                                                                                                                                                                                                                                                                                                                                                                                                                                                                                                                                                                                                                                                                                                                                                                                                                                                                                                                                                                                                                                                                                                                                                                                                                                                                                                                                                                                                                                                                                                                                                                                                                                                                                                                                                                                                                                                                                                                                                                                                                                                                                                                                                                                                                                                                                                                                                                                                                                                                                                                                                                                                                                                                                                                                                                                                                                                                                                                                                                                                                                                                                                                                                                                                                                                                                                                                                                                                                                                                                                                                                                                                                                                                           $ h_{1}'\cdots  h_{j-1}' h_{l_{1}}'\cdots h_{l_{\eta}}' h_{m}'\in\mathfrak{L}(V)$ by induction hypotheses of ${\rm (ii)}$.
If $\beta\geq1$,  then
$\tilde{p}_{ h_{j}' h_{g_{1}}'\cdots  h_{g_{\beta-1}}',  h_{g_{\beta}}'}\neq1$ and
$\tilde{p}_{ h_{1}'\cdots  h_{j-1}' h_{l_{1}}'\cdots h_{l_{\eta}}' h_{m}',  h_{g_{\beta}}'}$
$\neq1$ by $\mathbf{Assertion (k )}$ ${\rm (C_2)}$.
$ h_{1}'\cdots  h_{j-1}' h_{l_{1}}'\cdots  h_{l_{\eta}}' h_{m}' h_{j}' h_{g_{1}}'\cdots  h_{g_{\beta}-1}' h_{g_{\beta}}'\in\mathfrak{L}(V)$ by Lemma \ref {2.1}.
$ h_{1} h_{2}\cdots  h_{m}\sim  h_{1}'\cdots  h_{j-1}' h_{l_{1}}'\cdots h_{l_{\eta}}' h_{m}' h_{j}' h_{g_{1}}'\cdots  h_{g_{\beta}-1}' h_{g_{\beta}}'$. It is a contradiction. Then $\beta=0$.
$\eta=m-j-1$,  i.e. $\tilde{p}_{ h_{j}',  h_{r}'}=1$ for all $r\in\{j+1, \ldots, m\}$.
So $\tilde{p}_{ h_{i}'\cdots  h_{j-1}',  h_{r}'}=1$ for all $r\in\{j+1, \ldots, m\}$ by $\mathbf{Assertion (k )}$ ${\rm (C_4)}$.
Assume that there exists $\theta\in\{1, \ldots, i-1\}$ such that $\tilde{p}_{ h_{\theta}',  h_{j}'}\neq1$. If
$\tilde{p}_{ h_{j}',  h_{1}'\cdots  h_{\theta-1}' h_{\theta+1}'\cdots h_{i-1}' h_{i}'\cdots  h_{j-1}'  h_{j+1}'\cdots  h_{m}'}=1$,
then $\tilde{p}_{ h_{j}',  h_{1}'\cdots  h_{j-1}'  h_{j+1}'\cdots  h_{m}'}\neq1$,  $ h_{1}'$ $\cdots $ $ h_{j-1}'$ 
$  h_{j+1}'$ $\cdots $ $ h_{m}'\in$ $\mathfrak{L}(V)$  by induction hypotheses of ${\rm (iv)}$ and (\ref {e2.4.100}). 
$ h_{1}'\cdots  h_{j-1}'$ $ h_{j+1}'$ $\cdots$ $  h_{m}'$ $ h_{j}'$
$\in\mathfrak{L}(V)$ by \cite [Lemma 4.12]{WZZ15a}. $ h_{1} h_{2}\cdots  h_{m}\sim h_{1}'\cdots  h_{j-1}'  h_{j+1}'\cdots  h_{m}' h_{j}'$. It is a contradiction. If
$\tilde{p}_{ h_{j}',  h_{1}'\cdots  h_{\theta-1}' h_{\theta+1}'\cdots h_{i-1}' h_{i}'\cdots  h_{j-1}'  h_{j+1}'\cdots  h_{m}'}\neq1$,
$ h_{1}'$ $\cdots$ $ h_{\theta-1}'$ $ h_{\theta+1}'$ $\cdots$ $  h_{i-1}'$ $ h_{i}'$ $\cdots$ $ h_{j-1}' $ 
$ h_{j+1}'$ $\cdots $ $ h_{m}' $ $\in $ $\mathfrak{L}(V)$  by induction hypotheses {\rm (iv)} and (\ref {e2.4.100}). Then
$ h_{\theta}' h_{1}'\cdots  h_{\theta-1}' h_{\theta+1}'\cdots  h_{i-1}' h_{i}' $ $\cdots $ $ h_{j-1}'  h_{j+1}'\cdots $ $  h_{m}' h_{j}'$
 $\in\mathfrak{L}(V)$ by Lemma \ref {2.1}.
$ h_{1} h_{2}\cdots  h_{m}\sim  h_{\theta}' h_{1}'\cdots h_{\theta-1}' h_{\theta+1}'\cdots  h_{i-1}' h_{i}'\cdots  h_{j-1}'  h_{j+1}'\cdots  h_{m}' h_{j}'$. It is a contradiction. Then
$\tilde{p}_{ h_{\theta}',  h_{j}'}=1$ for all $\theta\in\{1, \ldots, i-1\}$,  which implies that ${\rm (C_2)}$ holds.

 ${\rm (C_3)}$ and ${\rm (C_4)}$ follows from ${\rm (C_2)}$.

($2^{\circ }$)  $A_{k+1} = A_k \cup \tau ^{k+1} (i)$ and $N_{k+1} = N_k - \tau ^{k+1} (i).$ If there exists $r\in A_{k+1} $ such that $\tilde{p}_{ h_{\alpha}',  h_{r}}= 1$ for all $\alpha\in\{i+1, \ldots, j\}$,  then
$\tilde{p}_{ h_{i}',   h_{r}}\neq 1$ by $\mathbf{Assertion (k )}$ and  $\tilde{p}_{ h_{i}'\cdots  h_{j}',   h_{r}}\neq 1$,  which  contradicts to $\mathbf{Assertion (k )}$. Thus for all $r\in A_{k+1} $, there exists $\alpha\in\{i+1, \ldots, j\}$ such that \begin {eqnarray} \label {e2.4.101}\tilde{p}_{ h_{\alpha}',   h_{r}}\neq 1 \end {eqnarray} and
$ h_{i+1}'\cdots  h_{j}' h_{j+1}'\cdots  h_{m}'\in\mathfrak L( V)$ by induction hypotheses of ${\rm (iv)}$.

Set $ \{g_1 , \cdots ,  g_{\beta},  l_1 ,  \cdots , l_{\eta}\} = \{\tau^{k+1}(2), \ldots, \tau^{k+1}(i-1)\}$  with $g_1 < \cdots < g_{\beta},  l_1 < \cdots < l_{\eta}$ such that
$\tilde{p}_{ h_{i}',   h_{g_{\lambda }}}\not= 1$ and $\tilde{p}_{ h_{i}',   h_{l_{\xi}}}= 1$ for all $ \lambda \in\{1, \ldots, \beta\}$,  all $\xi\in\{1, \ldots, \eta\}$.
Then $ h_{l_{1}}'\cdots  h_{l_{r}}' h_{i+1}'\cdots  h_{m}'\in\mathfrak L( V)$
for all $r\in\{1, \ldots, \eta\}$ and $ h_{g_{1}}'\cdots  h_{g_{t}}' h_{i}'\in\mathfrak L( V)$
for all $t\in\{1, \ldots, \beta\}$ by induction hypotheses of ${\rm (iv)}$.
If $\tilde{p}_{ h_{i}',   h_{1}'}\not= 1$,  then
$\tilde{p}_{ h_{l_{1}}'\cdots  h_{l_{\eta}}' h_{i+1}'\cdots  h_{m}',   h_{1}'}\not= 1$ and
$\tilde{p}_{ h_{g_{1}}'\cdots  h_{g_{\beta}}' h_{i}',   h_{1}'}\not= 1$ by $\mathbf{Assertion (k )}$ ${\rm (C_4)}$. Then $ h_{1}' h_{g_{1}}'\cdots  h_{g_{\beta}}' h_{i}' h_{l_{1}}'\cdots  h_{l_{\eta}}' h_{i+1}'\cdots  h_{m}'
\in  \mathfrak L( V)$ by Lemma \ref {2.1},
$ h_{1}\cdots  h_{m}\sim  h_{1}'\cdots  h_{m}'
\sim  h_{1}' h_{g_{1}}'\cdots  h_{g_{\beta}}' h_{i}'  h_{l_{1}}'\cdots h_{l_{\eta}}' h_{i+1}'\cdots  h_{m}'$.  It is a contradiction. Thus  $\tilde{p}_{ h_{i}',   h_{1}'}= 1$.
$ h_{1}' h_{l_{1}}'\cdots  h_{l_{\eta}}' h_{i+1}'\cdots h_{m}'\in\mathfrak{L}(V)$ by induction hypotheses of ${\rm (iv)}$.
If $\beta\geq1$,  then
$\tilde{p}_{ h_{g_{2}}'\cdots  h_{g_{\beta}}' h_{i}',  h_{g_{1}}'}\neq1$ and
$\tilde{p}_{ h_{1}' h_{l_{1}}'\cdots  h_{l_{\eta}}' h_{i+1}'\cdots h_{m}',  h_{g_{1}}'}$
$\neq1$ by $\mathbf{Assertion (k )}$ ${\rm (C_3)}$. then
$ h_{g_{1}}' h_{g_{2}}'\cdots  h_{g_{\beta}}' h_{i}' h_{1}' h_{l_{1}}'\cdots h_{l_{\eta}}' h_{i+1}'\cdots  h_{m}'\in\mathfrak{L}(V)$ by Lemma \ref {2.1}.
$ h_{1} h_{2}\cdots  h_{m}\sim  h_{g_{1}}' h_{g_{2}}'\cdots $ $ h_{g_{\beta}}' h_{i}' h_{1}' h_{l_{1}}'\cdots  h_{l_{\eta}}' h_{i+1}'\cdots  h_{m}'$, which   is a contradiction. Then $\beta=0$.
$\eta=i-2$,  i.e. $\tilde{p}_{ h_{i}',  h_{r}'}=1$ for all $r\in\{1, \ldots, i-1\}$.
So $\tilde{p}_{ h_{i+1}'\cdots  h_{j}',  h_{r}'}=1$ for all $r\in\{1, \ldots, i-1\}$ by $\mathbf{Assertion (k )}$ ${\rm (C_4)}$.
Assume that there exists $\theta\in\{j+1, \ldots, m\}$ such that $\tilde{p}_{ h_{\theta}',  h_{i}'}\neq1$. If
$\tilde{p}_{ h_{i}',  h_{1}'\cdots  h_{i-1}' h_{i+1}'\cdots  h_{j}' h_{j+1}'\cdots h_{\theta-1}' h_{\theta+1}'\cdots  h_{m}'}=1$,
then $\tilde{p}_{ h_{i}',  h_{1}'\cdots  h_{i-1}'  h_{i+1}'\cdots  h_{m}'}\neq1$   and $h_{1}'\cdots  h_{i-1}'  h_{i+1}'\cdots  h_{m}'
\in\mathfrak{L}(V)$  by induction hypotheses ${\rm (iv)}$ and (\ref {e2.4.101}).  $ h_{i}' h_{1}'\cdots h_{i-1}'  h_{i+1}'\cdots  h_{m}'
$ $\in \mathfrak{L}(V)$ by \cite [Lemma 4.12]{WZZ15a}. $ h_{1} h_{2}\cdots  h_{m}\sim h_{i}' h_{1}'\cdots  h_{i-1}'  h_{i+1}'\cdots  h_{m}'$, which  is a contradiction. 
\\
If $\tilde{p}_{ h_{i}',  h_{1}'\cdots  h_{i-1}'  h_{i+1}'\cdots  h_{j}' h_{j+1}'\cdots h_{\theta-1}' h_{\theta+1}'\cdots  h_{m}'}\neq1$, then
$ h_{1}'\cdots  h_{i-1}'  h_{i+1}'\cdots  h_{j}' h_{j+1}'\cdots h_{\theta-1}' h_{\theta+1}'$ $\cdots $ $ h_{m}'
\in\mathfrak{L}(V)$  by induction hypotheses  of {\rm (iv)} and  (\ref {e2.4.101}).
$ h_{i}' h_{1}'\cdots  h_{i-1}'  h_{i+1}'\cdots $ $ h_{j}' h_{j+1}' $ $\cdots $ $ h_{\theta-1}' h_{\theta+1}' $ $\cdots $ $ h_{m}' h_{\theta}'
$
 $\in\mathfrak{L}(V)$ by Lemma \ref {2.1}.
$ h_{1} h_{2}\cdots  h_{m}\sim  h_{i}' h_{1}'$ $\cdots $ $ h_{i-1}'$ $  h_{i+1}'$ $\cdots $ $h_{j}' h_{j+1}'\cdots  h_{\theta-1}' h_{\theta+1}'\cdots $ $ h_{m}' h_{\theta}'$, which  is a contradiction. Then
$\tilde{p}_{ h_{\theta}',  h_{i}'}=1$ for all $\theta\in\{j+1, \ldots, m\}$,
which implies  ${\rm (C_2)}$.

 ${\rm (C_3)}$ and ${\rm (C_4)}$ follows from ${\rm (C_2)}$.

Step 3. In  $\mathbf{Assertion (m-2)}$,  It is a contradiction by ${\rm (C_3)}$ and ${\rm (C_4)}$.

{\rm (iii)} By   Lemma \ref {2.2},  $u$ is connected. which implies that there exist two connected monomials $v$ and $w$
such that $u\sim  vw$. By inductive assumption,   $v$ and $w$ belong to $\mathfrak L(V).$
Consequently,  {\rm (iii)} holds by {\rm (i)}.

{\rm (iv)} By Lemma \ref {2.3},  $u \sim vw$ such that $v$ and $w$ are connected,
as well as,   $v\diamondsuit  w.$  By inductive assumption,
$v,  w \in  \mathfrak{L}(V). $ It follows from {\rm (i)} and {\rm (ii)}.  \hfill $\Box$

\begin {Corollary} \label {2.5} If $ h_{1} h_{2}\cdots  h_{m}\neq0$ and $p_{h_1,   h_1} \not= 1$ when   $\mu (h_{1} h_{2}\cdots  h_{m}) =h_1$,  then $ h_{1} h_{2}\cdots  h_{m}\in\mathfrak{L}(V)$
if and only if $\mu ( h_{1} h_{2}\cdots  h_{m})$ is connected.
\end {Corollary}
\noindent {\it Proof.}  It follows from Lemma \ref {2.2} and Theorem \ref {2.4} {\rm (iv)}.
\hfill $\Box$

\section { A basis of $\mathfrak L(V) $ }

\begin {Lemma}  \label{3.1} Assume that $V$ is  a  braided vector space   of diagonal type. If $u\neq0$ are homogeneous elements in $\mathfrak B(V)$ and $p_{i,   i} \not= 1$ when   $\mu (u) =x_i$
, $ i\in\{1,\ldots,n\}$,  then $\mu (u)$  is connected (i.e. every monomial of $u$ is connected ) if and only if $u\in \mathfrak L(V)$.
\end {Lemma}

\noindent {\it Proof.}  The necessity follows from  Corollary \ref {2.5}. We now prove the sufficiency. If $u$ is not connected and $u = \sum \limits _{i=1}^r k_i \sigma _i (u_i)$ with $k_i \in F^*$,  where  $u_i$ is a non-zero disconnected monomial and
 $\sigma _i$ is a method of bracket on letters of $u_i$ for $1\le i \le r$. By Lemma \ref {2.2},  $\sigma _i (u_i)=0$ for $1\le i \le r$ and $u=0$,  which is a contradiction.
 \hfill $\Box$


\begin {Lemma} \label {3.2} If $[u] \in D$,  then $u$ is connected and $u \not=0$.
\end {Lemma}
\noindent {\it Proof.}  By \cite [Cor. 1]{Kh99},  $u\not= 0.$ Obviously  $[u] \not=0$ and $[u]\in\mathfrak L( V)$. Considering Lemma \ref {3.1}  we complete the proof. \hfill $\Box$

\begin {Theorem} \label {3.3} If  $\mathfrak B(V) $ is Nichols algebra of diagonal type with $\dim V\geq2$ and $p_{i,   i} \not=  1$  for $1\le i\le n$,  then the set $\{[u_{1}]^{k_1}[u_{2}]^{k_2}\cdots  [u_{s}]^{k_s}\  \mid \  [u_{i}]\in D,  \mid D \mid = s
; 0 \le k_i <  h_{u_i};  1\le i \le s;   u_s<u_{s-1}<\cdots< u_1, \mu ([u_{1}]^{k_1}[u_{2}]^{k_2}\cdots  [u_{s}]^{k_s}) \hbox { is connect }, \sum \limits_{i=1}^{s}k_{i}>0\}$ is a basis of $\mathfrak L (V)$.
\end {Theorem}
\noindent {\it Proof.} It follows from \cite [Th. 1.4.6]{He05},   and Corollary \ref {2.5} and Lemma \ref {3.1}. \hfill $\Box$

\section {Dimension of $\mathfrak L(V) $}

In this section we give the basis for  $\mathfrak L(V) $ of
 arithmetic root systems and  obtain the dimensions of $\mathfrak L(V) $ of finite Cartan type.

Let $V_{i_1,  \cdots,  i_r}$ denote the braided vector subspace generated by $\{x_{i_1}, x_{i_2},  \cdots,   x_{i_r} \}$ of $V$; $D_{i_1,  \cdots,  i_r}$ denote $ \{[u] \mid [u] \hbox { is a hard super-letter of  } \mathfrak B(V_{i_1,  i_2,  \cdots,  i_r})   \}$;

\noindent $ L_{i_1,  \cdots,  i_r} := \{[u_{1}]^{k_1}[u_{2}]^{k_2}\cdots  [u_{s}]^{k_s}\  \mid \  [u_{j}]\in D_{i_1,  \cdots,  i_r}; 0 \le k_j < {\rm ord} (p_{u_j,  u_j});  1\le j \le s; $

\noindent $\mid D_{i_1,  \cdots,  i_r} \mid = s;  u_s<u_{s-1}<\cdots< u_1,  \mu ([u_{1}]^{k_1}[u_{2}]^{k_2}\cdots  [u_{s}]^{k_s}) \hbox { is connected }, \sum \limits_{j=1}^{s}k_{j}>0\}$.

\noindent $  B_{i_1,  \cdots,  i_r} := \{[u_{1}]^{k_1}[u_{2}]^{k_2}\cdots  [u_{s}]^{k_s}\  \mid \  [u_{j}]\in D_{i_1,  \cdots,  i_r}; 0 \le k_j < {\rm ord} (p_{u_j,  u_j});  1\le j \le s; $

\noindent $\mid D_{i_1,  \cdots,  i_r} \mid = s, \sum \limits_{j=1}^{s}k_{j}>0 \}$; Let $V_{s;  t}$ denote $V_{s,  s+1,   \cdots,  t}$
 in short; similarly we have  $B_{s;  t}$ and  $L_{s; t}$. Let $ B_{i;  j} := \emptyset $ and  $ L_{i;  j} := \emptyset $ when $i >j.$

\begin {Lemma} \label {4.1} Assumed that $\mathfrak B(V) $ is connected Nichols algebra of diagonal type with $\dim V>2$ and $\Delta(\mathfrak B(V)) $ is an arithmetic root system. The following hold.

{\rm (i)} ($n \ge 1$) If pure  generalized Dynkin graph  is

$
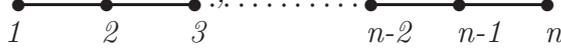
\begin{picture}(100,      15)
\put(27,      1){\makebox(0,     0)[t]{$\bullet$}}
\put(60,      1){\makebox(0,      0)[t]{$\bullet$}}
\put(93,     1){\makebox(0,     0)[t]{$\bullet$}}
\put(159,      1){\makebox(0,      0)[t]{$\bullet$}}
\put(192,     1){\makebox(0,      0)[t]{$\bullet$}}
\put(225,     1){\makebox(0,     0)[t]{$\bullet$}}
\put(28,      -1){\line(1,      0){30}}
\put(61,      -1){\line(1,      0){30}}
\put(130,     1){\makebox(0,     0)[t]{$\cdots\cdots\cdots\cdots$}}
\put(160,     -1){\line(1,      0){30}}
\put(193,      -1){\line(1,      0){30}}
\put(22,     -15){1}
\put(58,      -15){2}
\put(91,      -15){3}
\put(157,      -15){n-2}
\put(191,      -15){n-1}
\put(224,      -15){n}
\end{picture}$, \\
\\
then
\begin {eqnarray} \label {ep4.2.1} L_{1; n}&=&B_{1; n}-\cup _{i=1}^{n-2}L_{1; i}B_{i+2; n} \\
&=&B_{1; n}-\cup _{i=2}^{n-2}(L_{1; i}- L_{1;  i-1})B_{i+2; n} - L_{1; 1}B_{3; n}\\
&=& B_{1; n}-\cup _{i=1}^{n-3}L_{1; i}(B_{i+2; n}-B_{i+3; n})-L_{1; n-2}B_{n; n}
. \end {eqnarray}
\begin {eqnarray} \label {eppp4.2.2}
 \mid L_{1; n}\mid &=&\mid B_{1; n}\mid -\sum \limits _{i=1}^{n-2} \mid L_{1; i} \mid(\mid B_{i+2; n}\mid-\mid B_{i+3; n}\mid).
\end  {eqnarray}
\begin {eqnarray} \label {epp4.2.2} =\sum \limits _{j=0}^{{\rm int}(\frac{n-1}{2}) }(-1)^j u_j,
\end {eqnarray}
where
$u _0= \mid B_{1; n}\mid $ and
$ u_j = \sum \limits _{n_1=1}^{n -2}\sum \limits _{n_2=1}^{n_{1} -2 }\cdots \sum \limits _{n_j=1}^{n_{j-1} -2} \mid B_{1; n_j}\mid ( \mid B_{n_j+2; n_{j-1}}\mid
-\mid B_{n_j+3; n_{j-1}}\mid   ) \cdots ( \mid B_{n_2+2; n_1}\mid
-\mid B_{n_2+3; n_1}\mid   )( \mid B_{n_1+2; n}\mid
-\mid B_{n_1+3; n}\mid   )$ for $j>0$.

{\rm (ii)} ($n \ge 4$) If pure  generalized Dynkin graph  is
\\

$
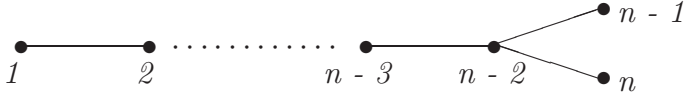
\begin{picture}(100,     15)
\put(27,      1){\makebox(0,     0)[t]{$\bullet$}}
\put(75,      1){\makebox(0,      0)[t]{$\bullet$}}
\put(156,     1){\makebox(0,     0)[t]{$\bullet$}}
\put(204,      1){\makebox(0,      0)[t]{$\bullet$}}
\put(245,     -11){\makebox(0,     0)[t]{$\bullet$}}
\put(245,    15){\makebox(0,     0)[t]{$\bullet$}}
\put(28,    -1){\line(1,    0){48}}
\put(115,    1){\makebox(0,     0)[t]{$\cdots\cdots\cdots\cdots$}}
\put(158,     -1){\line(1,     0){48}}
\put(202,    -1){\line(3,     1){42}}
\put(245,      -14){\line(-3,     1){42}}
\put(20,    - 15){1}
\put(70,    - 15){2}
\put(140,     - 15){n - 3}
\put(190,     - 15){n - 2}
\put(250,     8){n - 1}
\put(250,      -18){n}
\end{picture}$\\
\\
then
\begin {eqnarray} \label {ep4.2.3} L_{1; n} &=& B_{1; n} - B_{n-1;  n-1} B_{n;  n} -\cup _{i=1}^{n-3}L_{1; i}B_{i+2; n} \\
&=& B_{1; n}-\cup _{i=2}^{n-3}(L_{1; i} -L_{1;  i-1})B_{i+2; n} -  L_{1;  1}B_{3;  n} - B_{n-1;  n-1} B_{n;  n} \\
&=& B_{1; n}- L_{1;  n-3}B_{n-1;  n} - B_{n-1;  n-1} B_{n;  n}- \cup _{i=1}^{n-3}L_{1; i}(B_{i+2; n} - B_{i+3; n}).
\end {eqnarray}

 \begin {eqnarray} \label {eppp4.2.2'}  \mid L_{1; n}\mid &=&\mid B_{1; n}\mid -\mid B_{n-1; n-1}\mid\mid B_{n; n}\mid-\mid L_{1; n-3}\mid\mid B_{n-1; n}\mid
\nonumber \\
&&
  -\sum \limits _{i=1}^{n-3} \mid L_{1; i} \mid(\mid B_{i+2; n}\mid-\mid B_{i+3; n}\mid). \end  {eqnarray}
where $\mid L_{1; i} \mid $ is obtained by the formula (\ref {epp4.2.2}) when $1\le i \le n-3$.

{\rm (iii)} ($n \ge 6$) If pure  generalized Dynkin graph  is
\\
\\
\\
$
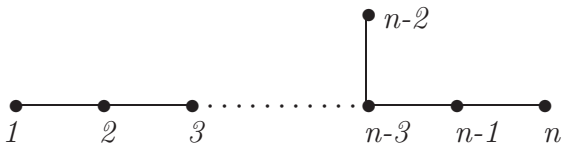
\begin{picture}(100,      15)
\put(27,      1){\makebox(0,     0)[t]{$\bullet$}}
\put(60,      1){\makebox(0,      0)[t]{$\bullet$}}
\put(93,     1){\makebox(0,     0)[t]{$\bullet$}}
\put(159,      1){\makebox(0,      0)[t]{$\bullet$}}
\put(192,     1){\makebox(0,      0)[t]{$\bullet$}}
\put(225,     1){\makebox(0,     0)[t]{$\bullet$}}
\put(28,      -1){\line(1,      0){30}}
\put(61,      -1){\line(1,      0){30}}
\put(130,     1){\makebox(0,     0)[t]{$\cdots\cdots\cdots\cdots$}}
\put(160,     -1){\line(1,      0){30}}
\put(193,      -1){\line(1,      0){30}}
\put(22,     -15){1}
\put(58,      -15){2}
\put(91,      -15){3}
\put(157,      -15){n-3}
\put(191,      -15){n-1}
\put(224,      -15){n}

\put(158,    1){\line(0,     1){33}}
\put(164,    28){n-2}
\put(155,    30){$\bullet$}

\end{picture}$\\

then  \begin {eqnarray} \label {ep4.2.5} L_{1; n}&=&B_{1; n}-\cup _{i=1}^{n-4}L_{1; i}B_{i+2; n} - L_{1; n-2}B_{n;  n} - B_{n-2; n-2}
 (B_{n-1; n}-B_{n;  n}) \\
 &=&B_{1; n}-\cup _{i=2}^{n-4}(L_{1; i}-L_{1;  i-1})B_{i+2; n} - L_{1; 1}B_{3;  n} \nonumber  \\
 && -
 (L_{1; n-2} -  L_{1; n-4})B_{n;  n} - B_{n-2; n-2}
 (B_{n-1; n}-B_{n;  n})
 . \end {eqnarray}
\begin {eqnarray} \label {ep4.2.6}
 \mid L_{1; n}\mid &=& \mid B_{1; n}\mid-\sum\limits_{i=2}^{n-4}(\mid  L_{1; i}\mid-\mid L_{1;  i-1}\mid) \mid B_{i+2; n}\mid   -\mid L_{1; 1}\mid \mid B_{3;  n} \mid \nonumber  \\
 && -
(\mid L_{1; n-2}\mid -  \mid L_{1; n-4}\mid)\mid B_{n;  n} \mid- \mid B_{n-2; n-2}\mid(\mid
 B_{n-1; n}\mid- \mid B_{n;  n} \mid). \end {eqnarray}
 where $\mid L_{1; i} \mid $ is obtained by the formula (\ref {epp4.2.2}) when $1\le i \le n-2$.

{\rm (iv)} ($n \ge 4$)  If pure  generalized Dynkin graph  is\\
\\
$
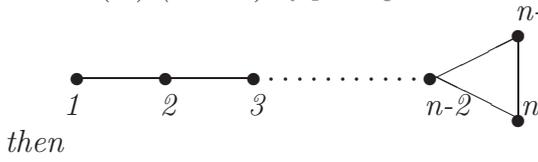
\begin{picture}(100,      15)
\put(27,      1){\makebox(0,     0)[t]{$\bullet$}}
\put(60,      1){\makebox(0,      0)[t]{$\bullet$}}
\put(93,     1){\makebox(0,     0)[t]{$\bullet$}}
\put(159,      1){\makebox(0,      0)[t]{$\bullet$}}
\put(192,     -15){\makebox(0,      0)[t]{$\bullet$}}
\put(192,    17){\makebox(0,      0)[t]{$\bullet$}}

\put(28,      -1){\line(1,      0){30}}
\put(61,      -1){\line(1,      0){30}}
\put(130,     1){\makebox(0,     0)[t]{$\cdots\cdots\cdots\cdots$}}
\put(192,     -16){\line(0,      1){30}}
\put(192,      -16){\line(-2,     1){33}}
\put(161,    -1){\line(2,     1){33}}

\put(22,     -15){1}
\put(58,      -15){2}
\put(91,      -15){3}
\put(157,      -15){n-2}
\put(193,      -15){n}
\put(191,    21){n-1}
\end{picture}$
\\
\\
then

\begin {eqnarray} \label {ep4.2.511} L_{1; n}&=&B_{1; n} - L_{1; 1} B_{3; n} - (L_{1; 2} - L_{1; 1}) B_{4; n} \cdots
- (L_{1; n-3} - L_{1; n-4})B_{n-1; n} \\
 &=&B_{1; n}-\cup _{i=2}^{n-3}(L_{1; i}-L_{1;  i-1})B_{i+2; n} - L_{1; 1}B_{3;  n}
 . \end {eqnarray}
\begin {eqnarray} \label {ep4.2.612}
 \mid L_{1; n}\mid &=& \mid B_{1; n}\mid-\sum\limits_{i=2}^{n-3}(\mid  L_{1; i}\mid-\mid L_{1;  i-1}\mid) \mid B_{i+2; n}\mid   -\mid L_{1; 1}\mid \mid B_{3;  n} \mid . \end {eqnarray}
 where $\mid L_{1; i} \mid $ is obtained by the formula (\ref {epp4.2.2}) when $1\le i \le n-2$.

{\rm (v)} ($n \ge 5$)  If pure  generalized Dynkin graph  is
\\

$
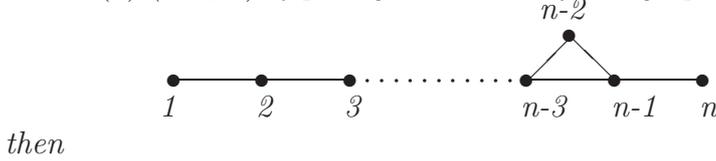
\begin{picture}(100,      15)
\put(27,      1){\makebox(0,     0)[t]{$\bullet$}}
\put(60,      1){\makebox(0,      0)[t]{$\bullet$}}
\put(93,     1){\makebox(0,     0)[t]{$\bullet$}}
\put(159,      1){\makebox(0,      0)[t]{$\bullet$}}
\put(192,     1){\makebox(0,      0)[t]{$\bullet$}}
\put(225,     1){\makebox(0,     0)[t]{$\bullet$}}
\put(175,    18){\makebox(0,     0)[t]{$\bullet$}}

\put(28,      -1){\line(1,      0){30}}
\put(61,      -1){\line(1,      0){30}}
\put(130,     1){\makebox(0,     0)[t]{$\cdots\cdots\cdots\cdots$}}
\put(160,     -1){\line(1,      0){30}}
\put(193,      -1){\line(1,      0){30}}
\put(191,      0){\line(-1,     1){17}}
\put(160,     -1){\line(1,      1){17}}

\put(22,     -15){1}
\put(58,      -15){2}
\put(91,      -15){3}
\put(157,      -15){n-3}
\put(191,      -15){n-1}
\put(224,      -15){n}
\put(164,    22){n-2}
\end{picture}$
\\
\\
then

\begin {eqnarray} \label {ep4.2.5121} L_{1; n}&=&B_{1; n} - L_{1; 1} B_{3; n} - (L_{1; 2} - L_{1; 1}) B_{4; n} \cdots
- (L_{1; n-5} - L_{1; n-6})B_{n-3; n}  \nonumber \\
 && - (L_{1; n-4} - L_{1; n-5})B_{n-2; n} - (L_{1; n-2} - L_{1; n-4})B_{n; n}\\
 &=&B_{1; n}-\cup _{i=2}^{n-4}(L_{1; i}-L_{1;  i-1})B_{i+2; n} - L_{1; 1}B_{3;  n}  -
 (L_{1; n-2} -  L_{1; n-4})B_{n;  n}
 . \end {eqnarray}
\begin {eqnarray} \label {ep4.2.6121}
 \mid L_{1; n}\mid &=& \mid B_{1; n}\mid-\sum\limits_{i=2}^{n-4}(\mid  L_{1; i}\mid-\mid L_{1;  i-1}\mid) \mid B_{i+2; n}\mid   -\mid L_{1; 1}\mid \mid B_{3;  n} \mid \nonumber  \\
 && -
(\mid L_{1; n-2}\mid -  \mid L_{1; n-4}\mid)\mid B_{n;  n} \mid. \end {eqnarray}
 where $\mid L_{1; i} \mid $ is obtained by the formula (\ref {epp4.2.2}) when $1\le i \le n-2$.

{\rm (vi)}  If pure  generalized Dynkin graph  is
\\
\\
$
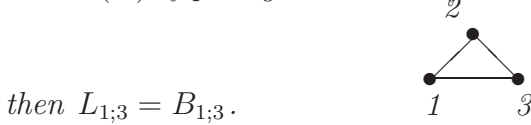
\begin{picture}(100,      15)
\put(159,      1){\makebox(0,      0)[t]{$\bullet$}}
\put(192,     1){\makebox(0,      0)[t]{$\bullet$}}
\put(175,    18){\makebox(0,     0)[t]{$\bullet$}}

\put(160,     -1){\line(1,      0){30}}

\put(191,      0){\line(-1,     1){17}}
\put(160,     -1){\line(1,      1){17}}

\put(157,      -15){1}
\put(191,      -15){3}

\put(164,    22){2}
\end{picture}$\\
 then $L_{1; 3}=B_{1; 3}$.

\end {Lemma}

\noindent {\it Proof.}
 {\rm (i)} We only determine which element in $B _{1;  n}$ is connected.
It is clear that the left hand of (\ref {ep4.2.1}) $\subseteq $  the right hand of (\ref {ep4.2.1}). If $u\in B _{1;  n}- L _{1;  n}$,   let $ i_u := {\rm min} \{ j \mid x_j \notin \mu (u) \hbox { and there exists } x_i \in \mu (u) \hbox { such that }  1\le i < j \le n\}$. By Lemma \ref {3.2},  there exist $v \in L_{1;  i_u-1}$ and  $w \in B_{i_u +1;  n }$ such that $u= vw$. Consequently,  the right hand of (\ref {ep4.2.1}) $\subseteq $  the left hand of (\ref {ep4.2.1}). therefore (\ref {ep4.2.1}) holds.

 \begin {eqnarray*}  \mid L_{1; n}\mid &=&\mid B_{1; n}\mid -\sum \limits _{n_1=1}^{n-2} \mid L_{1; n_1} \mid(\mid B_{n_1+2; n}\mid-\mid B_{n_1+3; n}\mid)\\
&=&  \mid B_{1; n}\mid \\
& & -\sum \limits _{n_1=1}^{n-2} \left(\mid B_{1; n_1}\mid -\sum \limits _{n_2=1}^{n_1-2} \mid L_{1; n_2} \mid(\mid B_{n_2+2; n_1}\mid-\mid B_{n_2+3; n_1}\mid)\right)\\
& &~~~~~~~~~\times \left(\mid B_{n_1+2; n}\mid-\mid B_{n_1+3; n}\mid\right)\\
&\cdots & \\
&=&\sum \limits _{j=0}^{{\rm int}(\frac{n-1}{2}) }(-1)^j u_j.
\end {eqnarray*}

Similarly,  we can show {\rm (ii)}-{\rm (vi)}.  \hfill $\Box$

\begin {Theorem} \label {4.2} The bases for  the Nichols braided Lie algebras of arithmetic root systems are given by Lemma \ref {4.1}.
\end {Theorem}
\noindent {\it Proof.} We can check that all pure generalized Dynkin diagrams in  \cite [Table A.1, A.2]{He05},  \cite [Table B, C]{He06a} are in Lemma \ref {4.1}. \hfill $\Box$

\vskip.1in
By \cite {Hu78},  $\mid D(A_n) \mid = C_{n+1} ^2, $ $\mid D(B_n) \mid = n^2 = \mid D(C_n) \mid, $ $\mid D(D_n) \mid = n^2-n,\mid D(E_6) \mid =36$,  $\mid D(E_7) \mid =63$,  $\mid D(E_8) \mid =120$,  $\mid D(F_4) \mid =24$,  $\mid D(G_2) \mid =6.$
By Lemma \ref {4.1} and \cite [Lemma 6.4] {WZZ15a}, we have the following results.

\begin {Theorem} \label {4.3}  Let ${\rm ord } (q) :=N$.\\

 {\rm (i)}
For $A_n$,  $n\ge 1$,
$\begin{picture}(100,      15)
\put(27,      1){\makebox(0,     0)[t]{$\bullet$}}
\put(60,      1){\makebox(0,      0)[t]{$\bullet$}}
\put(93,     1){\makebox(0,     0)[t]{$\bullet$}}
\put(159,      1){\makebox(0,      0)[t]{$\bullet$}}
\put(192,     1){\makebox(0,      0)[t]{$\bullet$}}
\put(28,      -1){\line(1,      0){30}}
\put(61,      -1){\line(1,      0){30}}
\put(130,     1){\makebox(0,     0)[t]{$\cdots\cdots\cdots\cdots$}}
\put(160,     -1){\line(1,      0){30}}
\put(22,     -15){1}
\put(58,      -15){2}
\put(91,      -15){3}
\put(157,      -15){n-1}
\put(191,      -15){n}
\put(22,     10){$q$}
\put(58,      10){$q$}
\put(91,      10){$q$}
\put(157,      10){$q$}
\put(191,      10){$q$}
\put(40,      5){$q^{-1}$}
\put(73,      5){$q^{-1}$}
\put(172,     5){$q^{-1}$}
\put(210,        -1)  {$, q \in F^{*}/\{1\}$. }
\end{picture}$\\
\\ then
$\dim \mathfrak L(V)=\sum \limits _{j=0}^{{\rm int}(\frac{n-1}{2})} (-1)^j u_j$,  where
$u _0= \mid B_{1;n}\mid $ and
$ u_j = \sum \limits _{n_1=1}^{n -2}\sum \limits _{n_2=1}^{n_{1} -2 }\cdots \sum \limits _{n_j=1}^{n_{j-1} -2} \mid B_{1; n_j}\mid ( \mid B_{n_j+2; n_{j-1}}\mid
-\mid B_{n_j+3; n_{j-1}}\mid   ) \cdots ( \mid B_{n_2+2; n_1}\mid
-\mid B_{n_2+3; n_1}\mid   )( \mid B_{n_1+2; n}\mid
-\mid B_{n_1+3; n}\mid   )$  for $j>0$
 and   $\mid B_{i;  k} \mid  = N ^{C_{k-i +2} ^2}-1$ for
$1\le i \le k \le  n$.
Furthermore,
\begin {eqnarray} \label {epppp4.2.2} && \dim \mathfrak L(V) = \mid L_{1; n} \mid  = N^ {C_{n+1}^2} -1
+
\sum \limits _{j=1}^{{\rm int}(\frac{n-1}{2})} (-1)^j \sum \limits _{n_1=1}^{n -2}\sum \limits _{n_2=1}^{n_{1} -2 }\cdots \sum \limits _{n_j=1}^{n_{j-1} -2}
( N ^{C_{n_j +1} ^2}-1  )\nonumber \\
&&
(N ^{C_{n_{j-1}- n_j } ^2} - N ^{C_{n_{j-1}- n_j -1} ^2})\cdots (N ^{C_{n_1- n_2 } ^2} - N ^{C_{n_1- n_2 -1} ^2})(N ^{C_{n- n_1 } ^2} - N ^{C_{n- n_1 -1} ^2}).
\end {eqnarray}
\\
{\rm (ii)} For $B_n$,  $n\ge 2$,
$\begin{picture}(100,      15)
\put(27,      1){\makebox(0,     0)[t]{$\bullet$}}
\put(60,      1){\makebox(0,      0)[t]{$\bullet$}}
\put(93,     1){\makebox(0,     0)[t]{$\bullet$}}
\put(159,      1){\makebox(0,      0)[t]{$\bullet$}}
\put(192,     1){\makebox(0,      0)[t]{$\bullet$}}
\put(225,     1){\makebox(0,     0)[t]{$\bullet$}}
\put(28,      -1){\line(1,      0){33}}
\put(61,      -1){\line(1,      0){30}}
\put(130,     -1){\makebox(0,     0)[t]{$\cdots\cdots\cdots\cdots$}}
\put(160,     -1){\line(1,      0){30}}
\put(193,      -1){\line(1,      0){30}}
\put(22,     -15){1}
\put(58,      -15){2}
\put(91,      -15){3}
\put(157,      -15){n-2}
\put(191,      -15){n-1}
\put(224,      -15){n}
\put(22,     10){$q^{2}$}
\put(58,      10){$q^{2}$}
\put(91,      10){$q^{2}$}
\put(157,      10){$q^{2}$}
\put(191,      10){$q^{2}$}
\put(224,      10){$q$}
\put(40,      5){$q^{-2}$}
\put(73,      5){$q^{-2}$}
\put(172,     5){$q^{-2}$}
\put(205,      5){$q^{-2}$}
\put(235,        -1)  {$, q \in F^{*}/\{1, -1\}$. }
\end{picture}$\\
\\  then
$\dim \mathfrak L(V)=\sum \limits _{j=0}^{{\rm int}(\frac{n-1}{2})} (-1)^j u_j$,  where
$u _0= \mid B_{1;n}\mid $ and
$ u_j = \sum \limits _{n_1=1}^{n -2}\sum \limits _{n_2=1}^{n_{1} -2 }\cdots \sum \limits _{n_j=1}^{n_{j-1} -2} \mid B_{1; n_j}\mid ( \mid B_{n_j+2; n_{j-1}}\mid
-\mid B_{n_j+3; n_{j-1}}\mid   ) \cdots ( \mid B_{n_2+2; n_1}\mid
-\mid B_{n_2+3; n_1}\mid   )( \mid B_{n_1+2; n}\mid
-\mid B_{n_1+3; n}\mid   )$ for $j>0$;   $\mid B_{i;  n} \mid  = N ^{(n-i +1)^2}-1$ for
$1\le i < n$ and $\mid B_{i;  k} \mid  = N ^{C_{k-i +2} ^2}-1$ for
$1\le i \le k <  n$,  when  N is  odd;
  $\mid B_{i;  n} \mid  = (\frac {N} {2}) ^{(n-i +1)^2-n +i-1}   N ^{n -i+1}-1$ for
$1\le i < n$ and $\mid B_{i;  k} \mid  = (\frac {N} {2}) ^{C_{k-i +2} ^2}-1$ for
$1\le i\le k < n$,  when  N is  even. Furthermore,
\begin {eqnarray} \label {epppp4.2.2'} && \dim \mathfrak L(V) = \mid L_{1; n} \mid  =N ^{n^2}-1
+
\sum \limits _{j=1}^{{\rm int}(\frac{n-1}{2})} (-1)^j \sum \limits _{n_1=1}^{n -2}\sum \limits _{n_2=1}^{n_{1} -2 }\cdots \sum \limits _{n_j=1}^{n_{j-1} -2}
( N ^{C_{n_j +1} ^2}-1  )\nonumber \\
&&
(N ^{C_{n_{j-1}- n_j } ^2} - N ^{C_{n_{j-1}- n_j -1} ^2})\cdots (N ^{C_{n_1- n_2 } ^2} - N ^{C_{n_1- n_2 -1} ^2}) \nonumber\\
 && (N ^{(n- n_1-1)^ 2 - n +n_1 +1 } - N ^{(n- n_1-2)^ 2 - n +n_1 +2 } ),
\end {eqnarray}
 when  N is  odd;
\noindent \begin {eqnarray} \label {epppp4.2.2'''} 
&& \dim \mathfrak L(V) = \mid L_{1; n} \mid  = (\frac {N} {2})^{n^2-n }  N ^n-1\\
& & 
+\sum \limits _{j=1}^{{\rm int}(\frac{n-1}{2})} (-1)^j \sum \limits _{n_1=1}^{n -2}\sum \limits _{n_2=1}^{n_{1} -2 }\cdots \sum \limits _{n_j=1}^{n_{j-1} -2}
( (\frac {N} {2}) ^{C_{n_j +1} ^2}-1  )\nonumber \\
&&
((\frac {N} {2}) ^{C_{n_{j-1}- n_j } ^2} - (\frac {N} {2}) ^{C_{n_{j-1}- n_j -1} ^2})\cdots ((\frac {N} {2}) ^{C_{n_1- n_2 } ^2} - (\frac {N} {2}) ^{C_{n_1- n_2 -1} ^2})\nonumber\\
 &&   ((\frac {N} {2}) ^{(n-n_{1} -1)^2-n +n_{1}+1}   N ^{n -n_{1}-1} - (\frac {N} {2}) ^{(n-n_{1} -2)^2-n +n_{1}+2}   N ^{n -n_{1}-2}),
\end {eqnarray}
 when  N is even.
\\ \\
{\rm (iii)} For $C_n$,  $n>2$,
$\begin{picture}(100,      15)
\put(27,      1){\makebox(0,     0)[t]{$\bullet$}}
\put(60,      1){\makebox(0,      0)[t]{$\bullet$}}
\put(93,     1){\makebox(0,     0)[t]{$\bullet$}}
\put(159,      1){\makebox(0,      0)[t]{$\bullet$}}
\put(192,     1){\makebox(0,      0)[t]{$\bullet$}}
\put(225,     1){\makebox(0,     0)[t]{$\bullet$}}
\put(28,      -1){\line(1,      0){33}}
\put(61,      -1){\line(1,      0){30}}
\put(130,     -1){\makebox(0,     0)[t]{$\cdots\cdots\cdots\cdots$}}
\put(160,     -1){\line(1,      0){30}}
\put(193,      -1){\line(1,      0){30}}
\put(22,     -15){1}
\put(58,      -15){2}
\put(91,      -15){3}
\put(157,      -15){n-2}
\put(191,      -15){n-1}
\put(224,      -15){n}
\put(22,     10){$q$}
\put(58,      10){$q$}
\put(91,      10){$q$}
\put(157,      10){$q$}
\put(191,      10){$q$}
\put(224,      10){$q^2$}
\put(40,      5){$q^{-1}$}
\put(73,      5){$q^{-1}$}
\put(172,     5){$q^{-1}$}
\put(205,      5){$q^{-2}$}
\put(235,        -1)  {$, q \in F^{*}/\{1, -1\}$. }
\end{picture}$\\
\\  then
$\dim \mathfrak L(V)=\sum \limits _{j=0}^{{\rm int}(\frac{n-1}{2})} (-1)^j u_j$,  where
$u _0= \mid B_{1;n}\mid  $ and
$ u_j = \sum \limits _{n_1=1}^{n -2}\sum \limits _{n_2=1}^{n_{1} -2 }\cdots \sum \limits _{n_j=1}^{n_{j-1} -2} \mid B_{1; n_j}\mid ( \mid B_{n_j+2; n_{j-1}}\mid
-\mid B_{n_j+3; n_{j-1}}\mid   ) \cdots ( \mid B_{n_2+2; n_1}\mid
-\mid B_{n_2+3; n_1}\mid   )( \mid B_{n_1+2; n}\mid
-\mid B_{n_1+3; n}\mid   )$ for $j>0$;   $\mid B_{i;  n} \mid  = N ^{(n-i +1)^2}-1$ for
$1\le i < n$ and $\mid B_{i;  k} \mid  = N ^{C_{k-i +2} ^2}-1$ for
$1\le i\le k < n$,  when  N is  odd;
  $\mid B_{i;  n} \mid  = N ^{(n-i +1)^2-n +i-1}   (\frac {N} {2}) ^{n -i +1}-1$ for
$1\le i < n$ and $\mid B_{i;  k} \mid  = N ^{C_{k-i +2} ^2}-1$ for
$1\le i\le k < n$,  when  N is  even.
Furthermore,
\begin {eqnarray} \label {epppp4.2.2''} && \dim \mathfrak L(V) = \mid L_{1; n} \mid  =N ^{n^2 }   -1
+
\sum \limits _{j=1}^{{\rm int}(\frac{n-1}{2})} (-1)^j \sum \limits _{n_1=1}^{n -2}\sum \limits _{n_2=1}^{n_{1} -2 }\cdots \sum \limits _{n_j=1}^{n_{j-1} -2}
( N ^{C_{n_j +1} ^2}-1  )\nonumber \\
&&
(N ^{C_{n_{j-1}- n_j } ^2} - N ^{C_{n_{j-1}- n_j -1} ^2})\cdots (N ^{C_{n_1- n_2 } ^2} - N ^{C_{n_1- n_2 -1} ^2})
\nonumber \\
&&(N ^{(n- n_1-1)^ 2 - n +n_1 +1 } - N ^{(n- n_1-2)^ 2 - n +n_1 +2 } ),
\end {eqnarray}
 when  N is  odd;
\begin {eqnarray} \label {epppp4.2.2''''} && \dim \mathfrak L(V) = \mid L_{1; n} \mid  =(\frac {N} {2}) ^{n }   N ^{n^2-n }-1\\
&&+
\sum \limits _{j=1}^{{\rm int}(\frac{n-1}{2})} (-1)^j \sum \limits _{n_1=1}^{n -2}\sum \limits _{n_2=1}^{n_{1} -2 }\cdots \sum \limits _{n_j=1}^{n_{j-1} -2}
( {N}  ^{C_{n_j +1} ^2}-1  )\nonumber \\
&&
(N ^{C_{n_{j-1}- n_j } ^2} - N ^{C_{n_{j-1}- n_j -1} ^2})\cdots (N ^{C_{n_1- n_2 } ^2} - N ^{C_{n_1- n_2 -1} ^2})\nonumber \\
&&  (N ^{(n-n_{1} -1)^2-n +n_{1}+1}   (\frac {N} {2}) ^{n -n_{1} -1} - N ^{(n-n_{1} -2)^2-n +n_{1}+2}   (\frac {N} {2}) ^{n -n_{1} -2}),
\end {eqnarray}
 when  N is even.
\\ \\
{\rm (iv)} For  $D_n$,  $n>3$,
$\begin{picture}(100,     15)
\put(12,      1){\makebox(0,     0)[t]{$\bullet$}}
\put(45,      1){\makebox(0,      0)[t]{$\bullet$}}
\put(111,     1){\makebox(0,     0)[t]{$\bullet$}}
\put(144,      1){\makebox(0,      0)[t]{$\bullet$}}
\put(170,     -11){\makebox(0,     0)[t]{$\bullet$}}
\put(170,    15){\makebox(0,     0)[t]{$\bullet$}}
\put(13,    -1){\line(1,    0){33}}
\put(80,    1){\makebox(0,     0)[t]{$\cdots\cdots\cdots\cdots$}}
\put(113,     -1){\line(1,     0){33}}
\put(142,    -1){\line(2,     1){27}}
\put(170,      -14){\line(-2,     1){27}}
\put(5,    - 15){1}
\put(40,    - 15){2}
\put(95,     - 15){n - 3}
\put(130,     - 15){n - 2}
\put(175,     8){n - 1}
\put(175,      -18){n}
\put(10,      10){$q$}
\put(45,      10){$q$}
\put(110,      10){$q$}
\put(143,      10){$q$}
\put(170,      -8){$q$}
\put(170,      18){$q$}
\put(25,      5){$q^{-1}$}
\put(123,      5){$q^{-1}$}
\put(150,     -4){$q^{-1}$}
\put(150,      10){$q^{-1}$}
\put(195,        -1)  {$, q \in F^{*}/\{1\}$. }
\end{picture}$\\
\\
then $ \dim \mathfrak L(V)=\mid B_{1; n}\mid -\mid B_{n-1; n-1}\mid\mid B_{n; n}\mid-\sum \limits _{i=1}^{n-3} \mid L_{1; i} \mid(\mid B_{i+2; n}\mid-\mid B_{i+3; n}\mid)-\mid L_{1; n-3}\mid\mid B_{n-1; n}\mid
$,
where
$\mid B_{ k ;n}\mid = N^{({n-k+1})^2-n +k -1}-1$ for $n-k +1 >2$,
$\mid B_{ n-1; n}\mid = N^3-1$, $ \mid B_{ n; n}\mid = \mid B_{ n-1; n-1}\mid = N-1$; $\mid L_{1; t} \mid $ is obtained by the formula  (\ref {epppp4.2.2}) when $1\le t \le n-3$.
Furthermore,
\begin {eqnarray} \label {epppp4.2.2'''''} && \dim \mathfrak L(V) = \mid L_{1; n} \mid  =
N^{n^2-n}-1 - (N-1)^2      -\sum \limits _{i=1}^{n-3} \left\{N^ {C_{i+1}^2} -1\right.\nonumber\\
&&+\sum \limits _{j=1}^{{\rm int}(\frac{i-1}{2})} (-1)^j \sum \limits _{n_1=1}^{n -2}\sum \limits _{n_2=1}^{n_{1} -2 } 
\cdots\sum \limits _{n_j=1}^{n_{j-1} -2}( N ^{C_{n_j +1} ^2}-1  )
(N ^{C_{n_{j-1}- n_j } ^2} - N ^{C_{n_{j-1}- n_j -1} ^2})\nonumber\\
&&\left.\cdots (N ^{C_{n_1- n_2 } ^2} - N ^{C_{n_1- n_2 -1} ^2})
(N ^{C_{i- n_1 } ^2} - N ^{C_{i- n_1 -1} ^2})\right\}\nonumber\\
&&(N^{({n-i-1})^2-n +i +1}-N^{({n-i-2})^2-n +i +2})\nonumber\\
&&-\left\{N^ {C_{n-2}^2} -1
+\sum \limits _{j=1}^{{\rm int}(\frac{n-4}{2})} (-1)^j \sum \limits _{n_1=1}^{n -5}\sum \limits _{n_2=1}^{n_{1} -2 }\cdots \sum \limits _{n_j=1}^{n_{j-1} -2}( N ^{C_{n_j +1} ^2}-1  )\right.\nonumber\\
&&(N ^{C_{n_{j-1}- n_j } ^2} - N ^{C_{n_{j-1}- n_j -1} ^2})
\cdots(N ^{C_{n_1- n_2 } ^2} - N ^{C_{n_1- n_2 -1} ^2})\nonumber\\
&&\left.(N ^{C_{n-3- n_1 } ^2} - N ^{C_{n-4- n_1 } ^2})\right\}(N^2-1).
\end {eqnarray}
\\
\\
{\rm (v)}  For $E_6$,
$\begin{picture}(100,    15)
\put(47,    1){\makebox(0,   0)[t]{$\bullet$}}
\put(87,    1){\makebox(0,0)[t]{$\bullet$}}
\put(127,   38){\makebox(0,0)[t]{$\bullet$}}
\put(127,    1){\makebox(0,    0)[t]{$\bullet$}}
\put(167,   1){\makebox(0,    0)[t]{$\bullet$}}
\put(207,   1){\makebox(0,   0)[t]{$\bullet$}}
\put(48,   -1){\line(1,    0){37}}
\put(88,  -1){\line(1,    0){37}}
\put(127,   1){\line(0,    1){37}}
\put(128,  -1){\line(1,    0){37}}
\put(168,  -1){\line(1,    0){37}}
\put(45,   -15){1}
\put(85,   -15){2}
\put(125,   -15){3}
\put(119,   26){4}
\put(166,   -15){5}
\put(206,   -15){6}
\put(49,    8){\makebox(0,   0)[t]{$q$}}
\put(89,    8){\makebox(0,0)[t]{$q$}}
\put(133,   40){\makebox(0,0)[t]{$q$}}
\put(129,    8){\makebox(0,    0)[t]{$q$}}
\put(169,   8){\makebox(0,    0)[t]{$q$}}
\put(209,   8){\makebox(0,   0)[t]{$q$}}
\put(60,   5){$q^{-1}$}
\put(101,   5){$q^{-1}$}
\put(128,   15){$q^{-1}$}
\put(143,   5){$q^{-1}$}
\put(180,   5){$q^{-1}$}
\put(227,       -1)  {$,q \in F^{*}/\{1\}$. }
\end{picture}$\\
\\then
\noindent \begin {eqnarray*}
  \dim \mathfrak L(V) &=& \mid B_{1; 6}\mid-\sum\limits_{i=2}^{6-4}(\mid  L_{1; i}\mid-\mid L_{1;  i-1}\mid) \mid B_{i+2; 6}\mid   -\mid L_{1; 1}\mid \mid B_{3;  6} \mid \nonumber  \\
 && -
(\mid L_{1; 6-2}\mid -  \mid L_{1; 6-4}\mid)\mid B_{6;  6} \mid- \mid B_{6-2; 6-2}\mid(\mid
 B_{6-1; 6}\mid- \mid B_{6;  6} \mid) \\
 &=&N ^{36}-1 - (\mid  L_{1;2}\mid-\mid L_{1; 1}\mid )  ( NN^3-1)   - (N-1) (N ^{C_5^2}-1) \\
 &-&
(\mid L_{1;4} \mid- \mid L_{1;2} \mid) (N-1)- (N-1) ( N ^3 - N ),
 \end {eqnarray*}
 where $\mid L_{1; t} \mid $ is obtained by the formula  (\ref {epppp4.2.2}) when $1\le t \le 4$.  \\
  \\
  \\
{\rm (vi)} For $E_7$,
$\begin{picture}(100,    15)
\put(37,    1){\makebox(0,    0)[t]{$\bullet$}}
\put(77,    1){\makebox(0,   0)[t]{$\bullet$}}
\put(117,    1){\makebox(0,0)[t]{$\bullet$}}
\put(157,   38){\makebox(0,0)[t]{$\bullet$}}
\put(157,    1){\makebox(0,    0)[t]{$\bullet$}}
\put(197,   1){\makebox(0,    0)[t]{$\bullet$}}
\put(237,   1){\makebox(0,   0)[t]{$\bullet$}}
\put(38,   -1){\line(1,    0){37}}
\put(78,   -1){\line(1,    0){37}}
\put(118,  -1){\line(1,    0){37}}
\put(157,   1){\line(0,    1){37}}
\put(158,  -1){\line(1,    0){37}}
\put(198,  -1){\line(1,    0){37}}
\put(35,   -15){1}
\put(75,   -15){2}
\put(115,   -15){3}
\put(155,   -15){4}
\put(149,   26){5}
\put(196,   -15){6}
\put(236,   -15){7}
\put(39,    8){\makebox(0,    0)[t]{$q$}}
\put(79,    8){\makebox(0,   0)[t]{$q$}}
\put(119,    8){\makebox(0,0)[t]{$q$}}
\put(163,   40){\makebox(0,0)[t]{$q$}}
\put(159,    8){\makebox(0,    0)[t]{$q$}}
\put(199,   8){\makebox(0,    0)[t]{$q$}}
\put(239,   8){\makebox(0,   0)[t]{$q$}}
\put(50,   5){$q^{-1}$}
\put(90,   5){$q^{-1}$}
\put(131,   5){$q^{-1}$}
\put(158,   15){$q^{-1}$}
\put(173,   5){$q^{-1}$}
\put(210,   5){$q^{-1}$}
\put(257,       -1)  {$,q \in F^{*}/\{1\}$. }
\end{picture}$\\
\\then
\noindent \begin {eqnarray*}
  \dim \mathfrak L(V) &=& \mid B_{1;7}\mid-\sum\limits_{i=2}^{7-4}(\mid  L_{1;i}\mid- \mid L_{1; i-1}\mid) \mid B_{i+2;7}\mid   -\mid L_{1;1}\mid \mid B_{3; 7} \mid \nonumber  \\
 && -
(\mid L_{1;7-2}\mid - \mid L_{1;7-4}\mid)\mid B_{7; 7} \mid- \mid B_{7-2;7-2}\mid
 (\mid B_{7-1;7}\mid - \mid B_{7; 7} \mid)\\
  &=& \mid B_{1;7}\mid-(\mid  L_{1;2}\mid-\mid L_{1; 1}\mid) \mid B_{3;7}\mid -(\mid  L_{1;3}\mid-\mid L_{1; 2}\mid) \mid B_{5;7}\mid\nonumber  \\
&&    -\mid L_{1;1}\mid \mid B_{3; 7} \mid 
  -(\mid L_{1;7-2}\mid - \mid L_{1;7-4}\mid)\mid B_{7; 7} \mid\\
&&- \mid B_{7-2;7-2}\mid (\mid B_{7-1;7}\mid- \mid B_{7; 7} \mid)\\
  &=& N^{63}-1
 -  ( \mid  L_{1;2}\mid-\mid L_{1; 1}\mid)  (N^{ 10}-1) -  (\mid  L_{1;3}\mid- \mid L_{1; 2}\mid)   (NN^{ 3}-1)
 \\
 &&   -
    ( N-1)( N^ {20}-1) -(\mid L_{1;5}\mid -\mid  L_{1;3} \mid) (N-1)-  (N-1)( N ^{3}
- N), \\
 \end {eqnarray*}
 where $\mid L_{1; t} \mid $ is obtained by the formula  (\ref {epppp4.2.2}) when $1\le t \le 5$.\\
 \\
\\
{\rm (vii)} For $E_8$,
$\begin{picture}(100,    15)
\put(27,    1){\makebox(0,    0)[t]{$\bullet$}}
\put(67,    1){\makebox(0,    0)[t]{$\bullet$}}
\put(107,    1){\makebox(0,   0)[t]{$\bullet$}}
\put(147,    1){\makebox(0,0)[t]{$\bullet$}}
\put(187,   38){\makebox(0,0)[t]{$\bullet$}}
\put(187,    1){\makebox(0,    0)[t]{$\bullet$}}
\put(227,   1){\makebox(0,    0)[t]{$\bullet$}}
\put(267,   1){\makebox(0,   0)[t]{$\bullet$}}
\put(28,   -1){\line(1,    0){37}}
\put(68,   -1){\line(1,    0){37}}
\put(108,   -1){\line(1,    0){37}}
\put(148,  -1){\line(1,    0){37}}
\put(187,   1){\line(0,    1){37}}
\put(188,  -1){\line(1,    0){37}}
\put(228,  -1){\line(1,    0){37}}
\put(22,  -15){1}
\put(65,   -15){2}
\put(105,   -15){3}
\put(145,   -15){4}
\put(185,   -15){5}
\put(179,   26){6}
\put(226,   -15){7}
\put(266,   -15){8}
\put(29,    8){\makebox(0,    0)[t]{$q$}}
\put(69,    8){\makebox(0,    0)[t]{$q$}}
\put(109,    8){\makebox(0,   0)[t]{$q$}}
\put(149,    8){\makebox(0,0)[t]{$q$}}
\put(193,   40){\makebox(0,0)[t]{$q$}}
\put(189,    8){\makebox(0,    0)[t]{$q$}}
\put(229,   8){\makebox(0,    0)[t]{$q$}}
\put(269,   8){\makebox(0,   0)[t]{$q$}}
\put(40,  5){$q^{-1}$}
\put(80,   5){$q^{-1}$}
\put(120,   5){$q^{-1}$}
\put(161,   5){$q^{-1}$}
\put(188,   15){$q^{-1}$}
\put(203,   5){$q^{-1}$}
\put(240,   5){$q^{-1}$}
\put(287,       -1)  {$,q \in F^{*}/\{1\}$. }
\end{picture}$\\
\\then
\begin {eqnarray*}
  \dim \mathfrak L(V) &=& \mid B_{1;8}\mid-\sum\limits_{i=2}^{8-4}(\mid  L_{1;i}\mid- \mid L_{1; i-1}\mid) \mid B_{i+2;8}\mid   -\mid L_{1;1}\mid\mid B_{3; 8} \mid \nonumber  \\
 && -
(\mid L_{1;8-2}\mid - \mid L_{1;8-4}\mid)\mid B_{8; 8} \mid- \mid B_{8-2;8-2}\mid(\mid
 B_{8-1;8}\mid-\mid B_{8; 8} \mid)\\
 &=&   \mid B_{1;8}\mid  -
(\mid  L_{1;2}\mid- \mid L_{1; 1}\mid )\mid B_{4;8}\mid
  -(\mid  L_{1;3}\mid-\mid L_{1; 2}\mid) \mid B_{5;8}\mid\\
   &&-
   (\mid  L_{1;4}\mid- \mid L_{1; 3}\mid )\mid B_{6;8}\mid
   -\mid L_{1;1}\mid \mid B_{3; 8} \mid\\
&&  -(\mid L_{1;6}\mid - \mid L_{1;4}\mid)\mid B_{8; 8} \mid
  -   \mid B_{6;6}\mid(\mid B_{7;8}\mid- \mid B_{8; 8}\mid)\\
 &=&  N  ^{120}-1  -
(\mid  L_{1;2}\mid- \mid L_{1; 1}\mid) ( N^ {20}-1)\\
&&  -(\mid  L_{1;3}\mid-\mid L_{1; 2}\mid )  ( N^ {C_5^2}-1)    \\
  &&-
  (\mid  L_{1;4}\mid-\mid L_{1; 3}\mid)  (N N ^3-1)
   - (N-1)(N^{36}-1)  \\
 && -
(\mid L_{1;6}\mid -\mid  L_{1;4} \mid )(N-1)-   (N-1)( N^3 - N )
 , \end {eqnarray*}
 where $\mid L_{1; t} \mid $ is obtained by the formula (\ref {epppp4.2.2}) when $1\le t \le 6$.\\
 \\
\\
{\rm (viii)} For $F_4$.
$\begin{picture}(100,      15)
\put(27,      1){\makebox(0,     0)[t]{$\bullet$}}
\put(60,      1){\makebox(0,      0)[t]{$\bullet$}}
\put(93,     1){\makebox(0,     0)[t]{$\bullet$}}
\put(126,      1){\makebox(0,    0)[t]{$\bullet$}}
\put(28,      -1){\line(1,      0){33}}
\put(61,      -1){\line(1,      0){30}}
\put(94,     -1){\line(1,      0){30}}
\put(22,     -15){1}
\put(58,      -15){2}
\put(91,      -15){3}
\put(124,     -15){4}
\put(22,     10){$q^2$}
\put(58,      10){$q^2$}
\put(91,      10){$q$}
\put(124,      10){$q$}
\put(40,      5){$q^{-2}$}
\put(73,      5){$q^{-2}$}
\put(106,     5){$q^{-1}$}
\put(145,        -1)  {$, q \in F^{*}/\{1, -1\}$. }
\end{picture}$\\
\\
 then
$\dim \mathfrak L(V)=\sum \limits _{j=0}^{1} (-1)^j u_j$,  where $u_0 =  \mid B_{1; 4}\mid $ and
$ u_1 = \sum \limits _{n_1=1}^{2} ( \mid B_{1; n_1}\mid ) ( \mid B_{n_1+2; 4}\mid
-\mid B_{n_1+3; 4}\mid   )$;  $\mid B_{1; 4}\mid  = N ^{24}-1$,  $\mid B_{1; 1}\mid = N -1$,
 $\mid B_{1; 2}\mid  = N ^{C_3^2}-1 = N^3-1$,  $\mid B_{3; 4 }\mid  = N ^{C_3^2}-1 = N^3-1$ and   $\mid B_{4; 4 }\mid  = N -1 $ when  $N$ is  odd;     $\mid B_{1; 4}\mid  = N ^{12} (\frac {N} {2}) ^{12}-1$,  $\mid B_{1; 1}\mid = \frac {N} {2} -1$,
 $\mid B_{1; 2}\mid  = (\frac {N} {2}) ^{C_3^2}-1 = (\frac {N} {2})^3-1$,  $\mid B_{3; 4 }\mid  = N ^{C_3^2}-1 = N^3-1$ and   $\mid B_{4; 4 }\mid  = N -1 $ when  $N$ is  even.
Furthermore
\begin {eqnarray} \label {epppp4.2.2''''''}  \dim \mathfrak L(V)   = N^ {24} -1 -( N -1)(N ^{3}-N ) - (N^3-1) (N-1),
\end {eqnarray} when $N$ is odd;
\begin {eqnarray} \label {epppp4.2.2.1}  \dim \mathfrak L(V)   = (\frac {N}{2})^ {12} N^{12} -1 -( \frac {N}{2} -1)(N ^{3}-N ) - ((\frac {N}{2})^ {3}-1) (N-1),
\end {eqnarray} when $N$ is even.
\\
\\
{\rm (ix)} For $G_2$,    $\begin{picture}(100,      15)
\put(27,      1){\makebox(0,     0)[t]{$\bullet$}}
\put(60,      1){\makebox(0,      0)[t]{$\bullet$}}
\put(28,      -1){\line(1,      0){33}}
\put(22,     -15){1}
\put(58,      -15){2}
\put(22,     10){$q$}
\put(58,      10){$q^3$}
\put(40,      5){$q^{-3}$}
  \ \ \ \ \ \ \ \ \ \ \ \ \ \ \ \ \ \ \ {$,q \in F^{*}/\{1, -1\}, q^3 \not= 1$. }
\end{picture}$\\
\\

 then
$\dim \mathfrak L(V)=N^ 6-1$ when $3 \nmid N$; $\dim \mathfrak L(V)=(\frac {N}{3})^3 N^3-1$ when $3 \mid N$.
 \end {Theorem}

\section { Non-zero monomials in  Nichols algebras }

Let $(s)_q := 1 + q + q^2 + \cdots + q ^{s-1}$ and $ (s) _q ! := (1)_q (2)_q \cdots (s)_q.$
\begin{Lemma}\label{5.1} Assume that $h_i \in \{x_1,  \cdots,  x_n\}$ for $1\le i \le m$ and $u _j$ is 1 or a monomial with $x_k \notin \mu (u_j)$ for $ 1\le j \le l+1$. Let $q:= p_{x_k, x_k} ^{-1}$ for convenience.

{\rm (i)}   \begin {eqnarray} \label {ep201} &&<y_{k}, u_{1}x_{k}u_{2}x_{k}\cdots u_{l}x_{k}u_{l+1}>
 \nonumber \\
&=& \sum \limits _{j=1}^{l}  q^{j-1} p _{k, u_1u_2\cdots u_j} ^{-1} u_1 x_k u_2 x_k \cdots x_k
(u_j u_{j+1}) x_k \cdots x_k u_{l+1}. \end {eqnarray}

{\rm (ii)} \begin {eqnarray} \label {ep202} <y_{k} ^l, u_{1}x_{k}u_{2}x_{k}\cdots u_{l}x_{k}u_{l+1}>
&=&  (l)_q!p_{k, u_{1}}^{-1}p_{k, u_{1}u_{2}}^{-1}\cdots p_{k, u_{1}u_{2}\cdots u_{l}}^{-1}
u_{1}u_{2}\cdots u_{l+1}
. \end {eqnarray}

{\rm (iii)} If ${\rm ord } (p_{h_i,h_i} ) >\mid {\rm deg}_ { h_i}( h_{1}\cdots  h_{m})\mid$ or
$p_{h_i,h_i} =1$ for all $i \in\{1, \ldots, m\}$,
then $ h_{1}\cdots  h_{m}\neq0$.
\end {Lemma}
\noindent {\it Proof.} {\rm (i)}  It can be obtained by induction on $l.$

 {\rm (ii)}
 We show this by induction on $l$. If $l=1$,  $<y_{k}, u_{1}x_{k}u_{2}>=p_{k, u_{1}}^{-1}
u_{1}u_{2}$. Assume $l>1.$  See that
\begin {eqnarray*}
 &&<y_{k}^{l}, u_{1}x_{k}u_{2}x_{k}\cdots u_{l}x_{k}u_{l+1}> \\
    &=&
<y_{k}^{l-1}, <y_{k}, u_{1}x_{k}u_{2}x_{k}\cdots u_{l}x_{k}u_{l+1}>> \\
&=&<y_k ^{l-1},  \sum \limits _{j=1}^{l}  q^{j-1} p _{k, u_1u_2\cdots u_j} ^{-1} u_1 x_k u_2 x_k \cdots x_k
(u_j u_{j+1}) x_k \cdots x_k u_{l+1}  > \\
&=&  \sum \limits _{j=1}^{l}  q^{j-1} p _{k, u_1u_2\cdots u_j} ^{-1}  < y_k ^{l-1},  u_1 x_k u_2 x_k \cdots x_k
(u_j u_{j+1}) x_k \cdots x_k u_{l+1}> \\
&=&  (l)_q!p_{k, u_{1}}^{-1}p_{k, u_{1}u_{2}}^{-1}\cdots p_{k, u_{1}u_{2}\cdots u_{l}}^{-1}
u_{1}u_{2}\cdots u_{l+1}.
\end {eqnarray*}

{\rm (iii)} We show this by induction on $t:= \mid \mu (h_1h_2\cdots h_m)\mid $. If $t=1$,  we obtain $ h_{1}\cdots  h_{m}\neq0$ by \cite [Lemma 1.3.3 {\rm (i)}]{He05}.
Assume $t >1$ and $ h_1h_2\cdots h_m = u_{1}x_{k}u_{2}x_{k}\cdots u_{l}x_{k}u_{l+1}$
 with $\mid \mu (u_{1}u_{2}\cdots u_{l+1})\mid $ $ = t-1$. Thus  $u_{1}u_{2}\cdots u_{l+1}\neq0$ by induction hypotheses,
$(1+p_{kk}^{-1})(1+p_{kk}^{-1}+p_{kk}^{-2})\cdots
(1+p_{kk}^{-1}+\cdots+p_{kk}^{-l+1})\neq0$ since ${\rm ord } (p_{x_{k}, x_{k}})>\mid {\rm deg}_ {x_k}( h_{1}\cdots  h_{m})\mid=l$ or $p_{h_i,h_i} =1$.
Hence $u_{1}x_{k}u_{2}x_{k}\cdots u_{l}x_{k}u_{l+1}\neq0$,  completing the proof.
 \hfill $\Box$

\begin {Corollary}  \label{5.2} If $u$ and $ v$ are monomials  in $\mathfrak B(V)$,  $\mu (u)\cap \mu (v)=\emptyset$,  $p_{ij} p_{ji}=1$ for any $x_{i} \in \mu (u)$,  $x_{j} \in \mu (v)$,  then the following conditions are equivalent:
{\rm (i)} $u\neq0, v\neq0$.
{\rm (ii)} $uv\neq0$.
{\rm (iii)} $uv\notin \mathfrak L(V)$.
\end {Corollary}
\noindent {\it Proof.} We know {\rm  (ii)} $\Longleftrightarrow$ {\rm  (iii)} by Lemma \ref {2.2}. {\rm (ii)}$\Longrightarrow$ {\rm (i)} is clear.
{\rm (ii)} $\Longleftarrow$ {\rm (i)}: Assume $\mid u\mid=k$, then $\exists y_{i_{1}},\ldots,y_{i_{k}}\in d(u)$ such that $<y_{i_{1}}\cdots y_{i_{k}},u>\in F^*$ and  $<y_{i_{1}}\cdots y_{i_{k}}, uv>$
$=<y_{i_{1}}\cdots y_{i_{k}}, u>v\neq0$.
 \hfill $\Box$

\begin {Lemma}  \label{5.3} If $u$ and $ v$ are two homogeneous elements  with   $\mu (u)\cap \mu (v)=\emptyset$, then  $uv =0$ implies $u=0$ or $v=0$.
\end {Lemma}
\noindent {\it Proof.} Without the lose of generality, there exists $i_0$ such that $\mu (v) \subseteq \{ 1, 2, \cdots i_0\}$ and $\mu (u) \subseteq \{ i_0+1, i_0+2, \cdots n\}$.
Assume $u \not=0$ and $v\not=0$ with $u = \sum \limits_{i=1} ^ s k_i u_i$  and $v = \sum \limits_{j=1} ^ t k_i' v_i$, where $k _i \not=0,$ $k _j' \not=0, $  $u_i  \in B_{1, i_0}$,  $v_j \in B_{i_0+1, n}$ for $1\le i \le i_0, $ $1+i_0\le j  \le n. $ Consequently, $uv = \sum \limits _{i=1} ^ s \sum \limits _{j=1} ^ t k_ik_j' u_iv_j \not= 0$ with $u_iv_j \in B_{1, n}$ and $k_ik_j \not=0$ for $1\le i \le s, $ $1\le j \le t. $  \hfill $\Box$

\section { Bases of Nichols Lie algebras $\mathfrak L^ - (V)$}

In this section we give the sufficient and necessary conditions for  $\mathfrak B(V) = F\oplus \mathfrak L^-(V)$ and $\mathfrak L^-(V)= \mathfrak L(V)$.
We also obtain an   explicit basis for $\mathfrak L^ - (V)$ over the quantum linear space $V$ with $\dim V=2$.

\vskip.1in
\subsection{Conditions for $\mathfrak B(V) = F\oplus \mathfrak L^-(V)$}

\begin{Lemma}\label{7.1} Assume that   $u_i$ is a homogeneous element and  $u_iu_j = r_{u_i, u_j} u_ju_i$  with $r_{u_i, u_j} \in F^*$ for $1\le i ,  j \le k$ and   $r_{u_i, u_j} =1$ when $u_i= u_j$.  Then
  \begin {eqnarray*}[u_{1 }, \cdots,  u_{m }]^{-}= (r_{u_{1 },u_{2 }\cdots u_{m }}-1)(r_{u_{2 },u_{3}\cdots u_{m }}-1)\cdots(r_{u_{m-1},u_{m }}-1)u_{m }\cdots u_{2 } u_1, \end  {eqnarray*} where $[u_{1 }, \cdots,  u_{m }]^{-}:= [u_1, [ u_2, \cdots [u_{m-1}, u_m]^- \cdots ]^-]^-$.

\end {Lemma}
\noindent {\it Proof.} Since $[u_{1 }, [u_{2 }, \cdots,  u_{m }]^{-}]=0$, we obtain $u_{1 }[u_{2 }, \cdots,  u_{m }]^{-}$ $=$ $r_{u_{1 }, u_{2 }\cdots u_{m }}[u_{2 }\cdots u_{m }]^{-}$ $u_{1 }$. Therefore,   \begin {eqnarray*}[u_{1 }, \cdots,  u_{m }]^{-} &=& u_{1 } [u_{2 }, \cdots,  u_{m }]^{-}-[u_{2 }, \cdots,  u_{m }]^{-}u_{1 }\\
&=& (r_{u_{1 },u_{2 }\cdots u_{m }}-1)[u_{2 }, \cdots,  u_{m }]^{-}u_{1 } \\
&=& (r_{u_{1 },u_{2 }\cdots u_{m }}-1)(r_{u_{2 },u_{3}\cdots u_{m }}-1)\cdots(r_{u_{m-1},u_{m }}-1)u_{m }u_{m-1}\cdots
u_{1 }. \ \ \hfill \Box  \end {eqnarray*}

 \begin {Lemma} \label {7.2}  Let $(L , [\ ]^-)$ be a Lie algebra and $u_1, u_2, \cdots, u_m \in L$.
 If $\sigma$ is a  method  of adding bracket $[\ ]^-$ on $u_1,  u_2,  \cdots,  u_m$,  then there exist some $\tau_{j}\in \mathbb S_{m}$, $\xi _{j}=1$ or $-1$ such that
 \begin {eqnarray} \label {ep4.1}  \sigma (u_1, u_2, \cdots, u_m)=\sum\limits_{j=1} ^r \xi _{j}[u_{\tau_{j}(1)}, \cdots,  u_{\tau_{j}(m)}]^{-}.\end  {eqnarray}

\end {Lemma}
 \noindent {\it Proof.}  We show  (\ref {ep4.1}) by induction on $m$. Obviously,  (\ref {ep4.1}) holds for $m =2$. Assume $m>2.$
 Let $\sigma (u_1, u_2, \cdots, u_m) = [\sigma_1 (u_1, u_2, \cdots, u_s), $ $\sigma _2(u_{s+1}, u_{s+2}, $ $\cdots,$ $ u_m)]^-$.

Now we show (\ref {ep4.1}) by induction on $s$.
 In case $s=1$, \begin   {eqnarray*} \sigma (u_1, u_2, \cdots, u_m) &=& [u_1, \sigma _2(u_{2}, u_{3}, \cdots, u_m)]^- \\
 &=&
 \sum\limits_{j=1} ^r \xi _{j} [ u_1, [u_{\tau_{j}(2)}\cdots u_{\tau_{j}(m)}]^{-} ]^-
\\
 &&(\hbox { \ by inductive assumption, where } \tau_j \in \mathbb S_{ \{2, 3. \cdots, m \}} \hbox  { for } 1\le j \le r)  \\
 &=& \sum\limits_{j=1} ^r \xi _{j} [ u_{\tau _j (1)}, u_{\tau_{j}(2)}\cdots u_{\tau_{j}(m)}]^{-}. \end  {eqnarray*} Therefore, (\ref {ep4.1}) holds. Assume $s>1$ and $\sigma_1 (u_1, u_2, \cdots, u_s) = [\sigma_3 (u_1, u_2, \cdots, u_k)$, $\sigma _4 (u_{k+1}, $ $u_{k+2}, $ $\cdots, u_s)]^-$.
 See  \\
 \begin   {eqnarray*} &&\sigma (u_1, u_2, \cdots, u_m) \\
  &=& [ [\sigma_3 (u_1, u_2, \cdots, u_k), \sigma _4 (u_{k+1}, u_{k+2}, \cdots, u_s)]^-  , \sigma _2(u_{s+1}, u_{s+2}, \cdots, u_m)]^- \\
 &=&[ [\sigma_3 (u_1, u_2, \cdots, u_k), \sigma _2 (u_{s+1}, u_{s+2}, \cdots, u_m)]^- , \sigma _4 (u_{k+1}, u_{k+2}, \cdots, u_s)]^-  \\
  && +[ \sigma_3 (u_1, u_2, \cdots, u_k), [ \sigma _4 (u_{k+1}, u_{k+2}, \cdots, u_s),  \sigma _2 (u_{s+1}, u_{s+2}, \cdots, u_m)]^- ]^- \\
 &=&
 \sum\limits_{j=1} ^{r_1} \xi _{j} [[ u_ {\tau _j (1)}, u_{\tau_{j}(2)}\cdots u_{\tau_{j}(k)}, u_ {\tau _j (s+1)}, u_{\tau_{j}(s+ 2)}\cdots u_{\tau_{j}(m)}]^-, \sigma _4 (u_{k+1}, u_{k+2}, \cdots, u_s)]^-  \\
 &&+
 \sum\limits_{j=r_1+1} ^{r_2} [ \sigma_3 (u_1, u_2, \cdots, u_k),  \xi _{j} [ u_ {\tau _j (k+1)}, u_{\tau_{j}(k+2)}\cdots u_{\tau_{j}(s)}, u_ {\tau _j (s+1)}, u_{\tau_{j}(s+ 2)}\cdots u_{\tau_{j}(m)}]^-]^-\\
  &&(\hbox { \ by first inductive assumption, where } \tau_j \in \mathbb S_{ \{2, 3. \cdots, k, s+1, s+2, \cdots, m \}} \hbox  { for } 1\le j \le r_1; \\
  &&
  \tau_j \in \mathbb S_{ \{k +1, k+2. \cdots,  m \}} \hbox  { for } r_1 +1\le j \le r_2)\\
  &=&\sum\limits_{j=r_2+ 1} ^{r_3} \xi _{j}[u_{\tau_{j}(1)}, \cdots,  u_{\tau_{j}(m)}]^{-}\\
  &&
  (\hbox { \ by second inductive assumption, where } \tau_j \in \mathbb S_{ \{1, 2. \cdots,  m \}} \hbox  { for } r_2+1\le j \le r_3). \end  {eqnarray*}
 Therefore, (\ref {ep4.1}) holds.  \hfill $\Box$

\begin {Proposition} \label {7.3} If $\mathfrak B(V) $ is a Nichols algebra of diagonal type, then $\mathfrak L(V)= \mathfrak L^-(V)$ if and only if $p_{ii}^2=1,p_{ij}=p_{ji}=1$ for $ 1\leq i\neq j\leq n$. In this case, $\mathfrak L(V)= \mathfrak L^-(V)=V$.
\end {Proposition}
\noindent {\it Proof.} The sufficiency. By Corollary \ref {2.5}, $\mathfrak L(V) = V$. By \cite  {WZZ15b}, $[x_i, x_j]^- =0$ for $i \not= j$ and  $\mathfrak L^-(V) = V$.

The necessity. We show this by following three steps.

{\rm (i)} $p_{ii}^2 =1$ since $x_i ^k \in \mathfrak L(V)$ for $k \le {\rm ord} (p_{ii})$ and $x_i ^m \notin \mathfrak L^-(V)$ when $0\not=x_i ^m$ and $m>1.$

{\rm (ii)}   If $p_{ij} p_{ji} =1$ and $p_{ij} \not= 1$ with $i \not= j,$ then $[x_i, x_j]^- \not=0$ and $[x_i, x_j]=0$. Consequently, $[x_i, x_j]^- \in \mathfrak L^-(V) - \mathfrak L(V), $ which is a contradiction.

{\rm (iii)} If $p_{ij} p_{ji} \not=1$ with $i < j$, then $0\not= [x_i, x_j] \in \mathfrak L(V)= \mathfrak L^-(V)$ and $[x_i, x_j] = k[x_i, x_j]^- $ with $0 \not= k \in F. $ By
Corollary \ref {2.5}, $x_jx_i \in \mathfrak L(V)$, which implies $x_jx_i = k'[x_i, x_j]^- $. This is a contradiction since $x_jx_i$ and $[x_i, x_j]$ are linearly  independent. $\Box$

\begin {Proposition} \label {7.4} Assume that $\mathfrak B(V) $ is a Nichols algebra of diagonal type. Then
$\mathfrak B(V) = F \oplus \mathfrak L^-(V)$
if and only if $p_{ii}=-1,p_{ij}p_{ji}=1$
for $\forall 1\leq i\neq j\leq n$
and  there  exist $\tau\in \mathbb  S_{m}$
such that $(p_{h_{\tau(1)},h_{\tau(2)}\cdots h_{\tau(m)}}-1)(p_{h_{\tau(2)},h_{\tau(3)}\cdots h_{\tau(m)}}-1)\cdots(p_{h_{\tau(m-1)},h_{\tau(m)}}-1)\neq0$
for $ \forall\ h_{1}>h_{2}>\cdots>h_{m}$
with $h_{i}\in \{x_{1},\ldots,x_{n}\},1\leq i\leq m$.
\end {Proposition}
\noindent {\it Proof.} The necessity.
If there exists $1\le i \le n$ such that  $ p_{ii}\neq-1$,
then $0\neq x_{i}^{2}\in \mathfrak B(V)$
and $x_{i}^{2}\notin \mathfrak L^-(V)$,
which  is a contradiction.
If there exist $i,  j $ such that $ p_{ij}p_{ji}\neq1$
with  $1\le i< j \le n$,
then $[x_{i},x_{j}]\neq0$ and $[x_{i},x_{j}]^-\neq0$.
Since $\mathfrak B(V) = F \oplus \mathfrak L^-(V)$, we have that there exist  $ k, k'\in F^{*}$ such that $[x_{i},x_{j}]=k[x_{i},x_{j}]^-$, and $x_{j}x_{i}=k'[x_{i},x_{j}]^-$,
which contradicts  to that $[x_{i},x_{j}]$ and
$x_{j}x_{i}$ are linearly independent.
Therefore, $V$ is a quantum linear space.

If there exist $ \ h_{1}>h_{2}>\cdots>h_{m}$
with $h_{i}\in \{x_{1},\ldots,x_{n}\},1\leq i\leq m$,  such that $(p_{h_{\tau(1)},h_{\tau(2)}\cdots h_{\tau(m)}}-1)(p_{h_{\tau(2)},h_{\tau(3)}\cdots h_{\tau(m)}}-1)\cdots(p_{h_{\tau(m-1)},h_{\tau(m)}}-1)= 0$ for any $\tau \in \mathbb S_m$.  By Lemma \ref {7.1} and Lemma \ref {3.2}, $\sigma ( h_{\tau (1)}, h_{\tau (2)}, \cdots, h_{\tau (m)} ) =0$ for any $\tau \in \mathbb S_m$ and  any  method $\sigma$ of adding bracket $[\ ]^-$ on $ h_{\tau (1)}, h_{\tau (2)}, \cdots, h_{\tau (m)}$. Consequently, $0\not= h_1h_2\cdots h_m \notin \mathfrak L^- (V)$, which is a contradiction.

The sufficiency. Obviously, $V$ is a quantum linear space. For $ \forall\ h_{1}>h_{2}>\cdots>h_{m}$
with $h_{i}\in \{x_{1},\ldots,x_{n}\},1\leq i\leq m$,
 there  exist $\tau\in \mathbb  S_{m}$
such that $(p_{h_{\tau(1)},h_{\tau(2)}\cdots h_{\tau(m)}}-1)(p_{h_{\tau(2)},h_{\tau(3)}\cdots h_{\tau(m)}}-1)\cdots(p_{h_{\tau(m-1)},h_{\tau(m)}}-1)\neq0$.
By Lemma \ref {7.1}, $h_1h_2\cdots h_m \sim h _{\tau (m)} h _{\tau (m-1)} $ $ \cdots $ $h _{\tau (1)} \in \mathfrak  L^-(V). $
 \hfill $\Box$

\subsection {Relationships between  graphs  and $\mathfrak L^ - (V)$}

Let
$\Gamma_a (V)_1 = \{ 1, 2, \cdots, n\}$ and $\Gamma_a (V)_2 = \{  a_{ij} \mid  p_{ij} \not= 1, i \not= j\}$. $\Gamma_a (V) = (\Gamma_a (V)_1, \Gamma_a (V)_2)$ is called the augmented Dynkin graph of $V$.

\begin {Proposition} \label {7.5} Assume that $h_i \in \{x_1,  \cdots,  x_n\}$ for $1\le i \le m$ and $u = h_1h_2\cdots h_m.$

{\rm (i)} If it is disconnected between monomial $u$ and monomial $v$ in $\Gamma_a (V)$ (i.e. ${p} _{x_i,  x_j} =1$  for any  $x_i \in \mu (u)$,  $x_j \in \mu (v)$ and $i\not= j$),  then $[u,  v]^- =0.$

{\rm (ii)} If $\Gamma_a (u)$ is  weakly disconnected  in $\Gamma_a (V)$ (i.e. $\Gamma_a (u)$ is disconnected in $\Gamma_a(V)$ or $ \mid \mu (u) \mid =1$  ),  then $\sigma (h_{1},  h_{2},  \cdots,   h_{m})=0$ for any  method $\sigma$ of adding bracket $[\ ]^-$ on $h_1,  h_2,  \cdots,  h_m.$

\end {Proposition}

\noindent {\it Proof.} {\rm (i)}  $u$ and $v$ are  commutative
(i.e. $uv= vu$ ) since $x_i x_j = x_jx_i$ for any  $x_i \in \mu (u)$,  $x_j \in \mu (v)$.

{\rm (ii)}
We show this by induction on $m$. $ [h_1,  h_2]^- =0$ for $m=2.$  For $m>2,  $ $\sigma (h_{1},  h_{2},  \cdots,  h_{m}) $

\noindent $= [\sigma _1 ( h_1 h_2 \cdots h_t ),  \sigma _2 ( h_{t+1} h_{t+2} \cdots h_m ) ]^-$.
If both $ h _1h_2\cdots h_t $ and $ h_{t+1} h_{t+2} \cdots h_m $ are connected  in $\Gamma_a (V)$,  then it is disconnected between  $h_1 h_2 \cdots h_t $ and $ h_{t+1} h_{t+2}$ $ \cdots$ $ h_m $  in $\Gamma_a (V)$. By Part {\rm (i)},   $\sigma (h_{1},  h_{2},  \cdots,  h_{m})=0$.
 If either  $h_1h_2 \cdots h_t $ or $ h_{t+1} h_{t+2} \cdots h_m$ is  disconnected  in $\Gamma_a (V)$,
then either  $\sigma _1 (h_1h_2\cdots h_t )=0$ or $ \sigma _2 (h_{t+1} h_{t+2} \cdots h_m)=0 $ by inductive hypothesis.

\hfill $\Box$

\subsection{ $\mathfrak L^ - (V)$ over quantum linear space $V$ with $\dim V=2$}\label {s7}

Let $R_m := \{ \alpha \mid \alpha$ is a primitive $m$-th root of $1 \}$.   For $a, b \in \mathbb Z$ with
$b \not=0$, if there exists $c\in \mathbb Z$ such that  $a= bc$,   then we say that $b$ is a  factor of $a$ and $a$ is divisible  by $b$, written $b \mid a$. Otherwise, we say that $a$ is not divisible  by $b$, written $b \nmid a$. For convenience, assume  $ b \nmid \infty $ and $ \infty \nmid b $.

\begin {Lemma} \label {7.6} Assume that   $u_i$ is a homogeneous element and  $u_iu_j = r_{u_i, u_j} u_ju_i$  with $r_{u_i, u_j} \in F^*$ for $1\le i ,  j \le k$ and   $r_{u_i, u_j} =1$ when $u_i= u_j$. Set  $l_{i}:=
\overline{l}_{u_{i}}$  and $r_{i,j}:=r_{u_{i},u_{j}}$ for $\forall\ 1\leq j, i\leq k$.

{\rm (i)}
\begin {eqnarray} \label {ep7.11.1} l_{k}^{m_{k}}l_{k-1}^{m_{k-1}}
\cdots l_{2}^{m_{2}}[u_{1}]^ - &=&(r_{2,1}-1)^{m_{2}}(r_{3,1}  r_{3, 2} ^{m_2}-1)^{m_{2}}\nonumber\\
&&\cdots(r_{k,1}r_{k,2}^{m_{2}}\cdots r_{k,k-1}^{m_{k-1}}-1)^{m_{k}}
u_{1}u_{2}^{m_{2}}\cdots u_{k}^{m_{k}} 
\end {eqnarray}
\noindent for $\forall\ m_{2}, \cdots, m_{k}\in\mathbb N$, i.e.
\begin {eqnarray} \label {ep7.11.2} l_{k}^{m_{k}}l_{k-1}^{m_{k-1}}
\cdots l_{2}^{m_{2}}[u_{1}]^ - =\lambda _2^{m_{2}}\lambda _3^{m_{3}}
\cdots\lambda _k^{m_{k}}
u_{1}u_{2}^{m_{2}}\cdots u_{k}^{m_{k}}
\end {eqnarray} where $\lambda _i :=  r_{i,1}r_{i,2}^{m_{2}}\cdots r_{i,i-1}^{m_{i-1}}-1$, $2\le i\le k$

{\rm (ii)} If $u_1, u_i \in \mathfrak L^- (V)$ and $\lambda _i \not=0$ for $2 \le i \le k, $ then $u_1u_2\cdots u_k \in \mathfrak L^- (V).$

{\rm (iii)} If $u_{2i} = u_2$ and $u_{2i-1} = u_1$ for $1\le i \le {\rm int} (\frac {k+1}{2})$, then
$ u_1u_2 ^ {m_1}\cdots u_k^ {m_k} \sim u_2 ^{\alpha _2} u_1^{\alpha_1}$;
$\lambda_k = r_{u_2, u_1} ^{\alpha _1} -1  $ with $ \alpha _1 = 1+ m_ 1 + m_2+ \cdots + m_{k-1}$  when $k$ is even;  $\lambda_k =  r_{u_1, u_2} ^{ \alpha_2} -1  $ with  $ \alpha _2 = m_ 2 + m_4+ \cdots + m_{k-1}$ when $k$ is odd.

{\rm (iv)} Assume that  $u_1, u_2 \in \mathfrak L^- (V)$ and $u_1^{\alpha_1}u_2^{\alpha_2}\not=0.$
${\rm ord } (r_{u_2, u_1}) \nmid  \alpha _1, $ or  ${\rm ord } (r_{u_1, u_2}) \nmid  \alpha _2, $
then
$u_1^{\alpha_1}u_2^{\alpha_2}\in \mathfrak L^- (V) $.

{\rm (v)} Assume that  $u_1, u_2 \in \{ x_1, x_2, \cdots, x_n\}$.  If  $u_1^{\alpha_1}u_2^{\alpha_2}\not=0,$ then
$u_1^{\alpha_1}u_2^{\alpha_2}\in \mathfrak L^- (V) $ if and only if ${\rm ord } (r_{u_2, u_1}) \nmid  \alpha _1, $ or  ${\rm ord } (r_{u_1, u_2}) \nmid  \alpha _2. $ \end {Lemma}

\noindent {\it Proof.} {\rm (i)} we show this  by induction on $k$.

For $k=1$, it is similar to \cite [Lemma 4.1 {\rm (iii)}] {WZZ15b}.

Assume that $k=2$. We show this by induction on $m_{2}$.  It is clear when $m_{2}=1$.
If $m_{2}>1$, then \begin {eqnarray*} l_{2}^{m_{2}}[u_{1}]^ - &=&[u_{2},l_{2}^{m_{2}-1}[u_{1}]^ -]^{-}
=(r_{21}-1)^{m_{2}-1}
[u_{2},u_{1}u_{2}^{m_{2}-1}]^{-}\\
&=&
(r_{21}-1)^{m_{2}-1}(r_{21}u_{1}u_{2}^{m_{2}}-u_{1}u_{2}^{m_{2}})=(r_{21}-1)^{m_{2}}
u_{1}u_{2}^{m_{2}}.\end  {eqnarray*}
Now  $k>2$.  We show this by induction on $m_{k}$. If $m_{k}=1$, then \begin {eqnarray*}
&& l_{k}l_{k-1}^{m_{k-1}}
\cdots l_{2}^{m_{2}}[u_{1}]^ -
=[u_{k},l_{k-1}^{m_{k-1}}
\cdots l_{2}^{m_{2}}[u_{1}]^ -]^ - \\
&=&(r_{2,1}-1)^{m_{2}}
(r_{3,1}r_{3,2}^{m_{2}}-1)^{m_{3}}\cdots(r_{k-1,1}r_{k-1,2}^{m_{2}}\cdots r_{k-1,k-2}^{m_{k-2}}-1)^{m_{k-1}}\\
&&[u_{k},u_{1}u_{2}^{m_{2}}\cdots u_{k-1}^{m_{k-1}} ]^{-} \\
&=&(r_{2,1}-1)^{m_{2}}
(r_{3,1}r_{3,2}^{m_{2}}-1)^{m_{3}}\cdots(r_{k-1,1}r_{k-1,2}^{m_{2}}\cdots r_{k-1,k-2}^{m_{k-2}}-1)^{m_{k-1}}\\
&& (r_{k,1}r_{k,2}^{m_{2}}\cdots r_{k,k-1}^{m_{k-1}}-1)u_{1}u_{2}^{m_{2}}\cdots u_{k-1}^{m_{k-1}}u_{k}.\end {eqnarray*}

If $m_{k}>1,$ then
\begin {eqnarray*} &&l_{k}^{m_{k}}l_{k-1}^{m_{k-1}}
\cdots l_{2}^{m_{2}}[u_{1}]^ -
=[u_{k},l_{k}^{m_{k}-1}l_{k-1}^{m_{k-1}}\cdots l_{2}^{m_{2}}[u_{1}]^ -]^ - \\
&&~~~=(r_{2,1}-1)^{m_{2}}
(r_{3,1}r_{3,2}^{m_{2}}-1)^{m_{3}}\cdots(r_{k-1,1}r_{k-1,2}^{m_{2}}\cdots r_{k-1,k-2}^{m_{k-2}}-1)^{m_{k-1}}\\
&&~~~ (r_{k,1}r_{k,2}^{m_{2}}\cdots r_{k,k-1}^{m_{k-1}}-1)
^{m_{k}-1} 
[u_{k},u_{1}u_{2}^{m_{2}}\cdots u_{k-1}^{m_{k-1}}u_{k}^{m_{k}-1} ]^{-}\\
&&~~~=(r_{2,1}-1)^{m_{2}}
(r_{3,1}r_{3,2}^{m_{2}}-1)^{m_{3}}\cdots(r_{k-1,1}r_{k-1,2}^{m_{2}}\cdots r_{k-1,k-2}^{m_{k-2}}-1)^{m_{k-1}}\\
&&~~~(r_{k,1}r_{k,2}^{m_{2}}\cdots r_{k,k-1}^{m_{k-1}}-1)^{m_{k}}
u_{1}u_{2}^{m_{2}}\cdots u_{k-1}^{m_{k-1}}u_{k}^{m_{k}}\\
&&~~  \hbox { \ \ (by induction hypotheses )}.
\end {eqnarray*}  Consequently, (\ref {ep7.11.1}) holds.

{\rm (ii)} It follows from {\rm (i)}.

{\rm (iii)}
It follows from {\rm (i)}.

{\rm (iv)}    If ${\rm ord } (r_{u_2, u_1}) \nmid  \alpha _1, $ then $\lambda _2 = (r_{u_2, u_1} ^{\alpha_1} -1) \not=0$ by Part {\rm (iii)} for $k=2$. Consequently, $u_1^{\alpha_1}u_2^{\alpha_2}\in \mathfrak L^- (V) $ by formula (\ref {ep7.11.2}).
Similarly, $u_1^{\alpha_1}u_2^{\alpha_2}\in \mathfrak L^- (V) $ when ${\rm ord } (r_{u_1, u_2}) \nmid  \alpha _2. $

{\rm (v)}  We only need show
the necessity. Assume that it does not holds. By {\rm (iii)}, $l_{k}^{m_{k}}l_{k-1}^{m_{k-1}}
\cdots l_{2}^{m_{2}}[u_{1}]^ -=0$ and $l_{k}^{m_{k}}l_{k-1}^{m_{k-1}}
\cdots l_{2}^{m_{2}}[u_{2}]^ -=0$ for any  $m_2, m_3, \cdots, m_k \in \mathbb N$. Considering Lemma \ref {3.2},  $u_1^{\alpha_1}u_2^{\alpha_2}=0$, which is a contradiction.
 \hfill $\Box$

\begin {Theorem} \label {7.7}  If  $\mathfrak L^ - (V)$ is  Nichols Lie algebra over quantum linear space $V$  with $\dim V =2$,
   then
    \begin {eqnarray} \label {e7.12.1}
    \mathfrak L^-(V)&=&{\rm span }\{x_{1},x_{2},x_{2}^{\alpha _{2}}x_{1}^{\alpha_{1}}\mid 1\leq \alpha_{2}< N_2
,1\leq \alpha _{1}< N_1,\nonumber\\
& &~~~~~~~~~~~ {\rm ord } (p_{12})\nmid  \alpha _{2} \hbox{ or }  {\rm ord } (p_{12})\nmid \alpha _{1} \},  
\end  {eqnarray} where  $N_i :=\left \{  \begin{array}{ll}
{\rm ord } (p_{ii}),  & \mbox {when }    {\rm ord } (p_{ii})>1 \\
\infty,   & \mbox {when }  {\rm ord } (p_{ii})=1  \\
\end{array}\right. $, for $1\le i \le 2.$
\end {Theorem}
\noindent {\it Proof.}
 Let $u_i \in \{x_1, x_2\}$   and  $r_{u_i, u_j} :=\left \{  \begin{array}{ll}
1,  & \mbox {when }    u_j= u_i \\
p_{u_i, u_j},   & \mbox {when }  u_j \not=u_i  \\
\end{array}\right. $. Then  $r_{u_i, u_j}r_{u_j, u_i}=1$ and $u_{i}u_{j}=r_{u_{i},u_{j}}u_{j}u_{i}$
for $1 \le i, j \le k$. Considering Lemma \ref {7.5} {\rm (v)} we complete the proof.  \hfill $\Box$

\begin {Example} \label {7.8} Assume $p_{1  2} = -1$ and    $ord(p_{1  1}) =3, ord(p_{2  2}) =5 $. Then  $ \mathfrak L^-(V)=F\{x_{1},x_{2},x_{2}x_{1},
x_{2}x_{1}^{2},x_{2}^{2}x_{1},x_{2}^{3}x_{1},
x_{2}^{3}x_{1}^{2},x_{2}^{4}x_{1}\}$. In this case, $ \mathfrak B (V)  =  \mathfrak L^ -  (V) \oplus span \{ x_i^ {m_i}  \mid  1< m_i \in \mathbb N \} \oplus F\oplus span \{x_{2}^{2}x_{1}^{2},x_{2}^{4}x_{1}^{2}\}.$
 \end {Example}

\section*{Acknowledgment}
 YZZ was partially supported by the Australian Research Council
through Discovery-Projects grant DP140101492.

\end {document}